\documentclass[9pt,shortpaper,twoside,web]{ieeecolor}
\usepackage{generic}
\usepackage{cite}
\usepackage{amsmath,amssymb,amsfonts}
\usepackage{algorithmic}
\usepackage{graphicx}
\usepackage{textcomp}
\usepackage{gensymb}
\usepackage{color}         
\usepackage{amsmath} 
\usepackage{amssymb}  
\usepackage{mathtools}
\mathtoolsset{showonlyrefs}
\usepackage[linesnumbered,ruled,vlined,algo2e]{algorithm2e}
\usepackage{algorithm}
\usepackage{comment}
\usepackage{ifthen}
\usepackage{amsmath}
\usepackage[bb=dsserif]{mathalpha} 
\usepackage{wrapfig}
\usepackage{lipsum} 
\usepackage[export]{adjustbox} 

\DeclareMathOperator*{\argmin}{arg\,min}
\DeclarePairedDelimiter{\norm}{\lVert}{\rVert}

\def\BibTeX{{\rm B\kern-.05em{\sc i\kern-.025em b}\kern-.08em
    T\kern-.1667em\lower.7ex\hbox{E}\kern-.125emX}}

\begin{document}
\newcommand{\gb}[1]{\textcolor{blue}{#1}}
\newcommand{\mh}[1]{\textcolor{red}{#1}}
\newcommand{\ml}[1]{\textcolor{green}{#1}}

\newcommand{\T}{\mathcal{T}}
\newcommand{\Y}{\mathcal{Y}}
\newcommand{\X}{\mathcal{X}}
\newcommand{\R}{\mathbb{R}}
\newcommand{\N}{\mathbb{N}}
\newcommand{\C}{\mathcal{C}}
\newcommand{\M}{\mathcal{M}}
\newcommand{\K}{\mathcal{K}}

\newcommand{\at}[2]{\alpha \left( {#1} ; t_{#2} \right)}
\newcommand{\bt}[2]{\beta \left( {#1} \right)}
\newcommand{\dt}[2]{\delta \left( {#1} \right)}
\newcommand{\pt}[2]{\rho^{#1}_{#2}}
\newcommand{\diam}{\textnormal{diam}}

\newcommand{\s}[3]{
\ifthenelse{\equal{#3}{}}{s(#2;\ind{#3})}{s_{#1}(#2)}
}
\newcommand{\q}[2]{
\ifthenelse{\equal{#2}{}}{q(#2)}{q_{#1}(#2)}
}
\newcommand{\Lx}{L_{x,\ell}}
\newcommand{\Ly}{L_{y,\ell}}
\newcommand{\Lj}{L_{J,\ell}}
\newcommand{\Lt}{L_{\ell}}

\newcommand{\f}[2]{f \big(#1;\ind{#2}\big) }
\newcommand{\g}[2]{g \big(#1;\ind{#2}\big) }
\newcommand{\J}[3]{J \big(#1,#2;\ind{#3}\big)}
\newcommand{\ind}[1]{t_{#1}}
\newcommand{\solSeries}{\left\{\x{*}{}{\ind{\ell}} \right\}_{\ind{\ell} \in \T}}
\newcommand{\Ci}{C_i}

\newcommand{\x}[3]{
\ifthenelse{\equal{#3}{}}{x^{#1}_{#2}}{x^{#1}_{#2}(#3)}
}
\newcommand{\y}[3]{
\ifthenelse{\equal{#3}{}}{y^{#1}_{#2}}{y^{#1}_{#2}(#3)}
}

\newtheorem{dfn}{Definition}
\newtheorem{theorem}{Theorem}
\newtheorem{assumption}{Assumption}
\newtheorem{lemma}{Lemma}
\newtheorem{problem}{Problem}
\newtheorem{remark}{Remark}
\newtheorem{corollary}{Corollary}

\title{Distributed Asynchronous Discrete-Time Feedback Optimization
}\author{Gabriel Behrendt, Matthew Longmire, Zachary I. Bell,
Matthew Hale, \IEEEmembership{Member, IEEE}
\thanks{
Gabriel Behrendt and Matthew Hale were supported by 
AFRL under grant FA8651-23-F-A006, 
AFOSR under grant FA9550-19-1-0169,  
and ONR under grant N00014-21-1-2495.
Matthew Longmire was supported by AFRL under contract FA8651-22-F-1052 under task orders FA8651-19-D-0037 and FAA8651-22-F-1045.
}
\thanks{Gabriel Behrendt, Matthew Longmire, and Matthew Hale are
with the Department of Mechanical and Aerospace Engineering
at the University of Florida, Gainesville, FL, USA.
Emails: \texttt{\{gbehrendt,m.longmire,matthewhale\}@ufl.edu} Zachary I. Bell is with AFRL/RW at Eglin AFB. Email: \texttt{zachary.bell.10@us.af.mil}
}
}

\maketitle

\begin{abstract}
In this article, we present an algorithm that drives the outputs of a network of agents to 
jointly track the solutions of time-varying optimization problems in a way that is robust 
to asynchrony in the agents' operations. We consider three operations that can be asynchronous: 
(1) computations of control inputs, (2) measurements of network outputs, and (3) communications 
of agents' inputs and outputs. We first show that our algorithm 
converges to the solution of a time-invariant feedback optimization problem in linear time. 
Next, we show that our algorithm drives outputs to track the solution of time-varying feedback optimization 
problems within a bounded error dependent upon the movement of the minimizers and degree of 
asynchrony in a way that we make precise. These convergence results are extended to quantify 
agents' asymptotic behavior as the length of their time horizon approaches infinity. Then, to 
ensure satisfactory network performance, we specify the timing of agents' operations relative 
to changes in the objective function that ensure a desired error bound. 
Numerical experiments confirm these developments and show the success
of our distributed feedback optimization algorithm under asynchrony. 
\end{abstract}

\section{Introduction}
\label{sec:introduction}
Time-varying optimization problems arise in machine learning, robotics, power systems, and others \cite{yang2009linear,koppel2015target,tang2019time}. These problems can model time-varying demands in power distribution systems~\cite{tang2017real} and robot navigation in cluttered dynamic environments~\cite{arslan2016exact}, among other engineering
problems. 
Time-varying optimization problems have been studied in both continuous-time~\cite{yong1993training,sun2017distributed,fazlyab2018prediction} and discrete-time~\cite{popkov2005gradient,tang2022running,simonetto2017prediction,bastianello2019prediction}, and methods for tracking their solutions include correction-only methods~\cite{colombino2020online,popkov2005gradient} and prediction-correction methods~\cite{simonetto2016class,bauman2004newton}. 
For a survey of time-varying optimization  see~\cite{simonetto2020time,hauswirth2021optimization}.

Time-varying optimization problems have been combined with control by embedding optimization algorithms
into feedback loops. This setup often uses the measured 
output of a dynamical system as the input to an optimization
algorithm. Then the optimization algorithm computes new control inputs for the system that
drive its outputs to track the time-varying solution of a time-varying optimization problem. 
The actual measured outputs of a system can be subject to disturbances, e.g., measurement
noise, and the use of measured outputs in this setup can provide robustness to such disturbances
without needing to explicitly estimate or model those disturbances~\cite{colombino2019towards}. 
In discrete time, measuring an output leads to a new optimization problem whose solution
is the optimal input at the next timestep. 
In some cases, the calculation of optimal inputs 
cannot be run to completion
due to practical constraints, 
e.g., a low-power computer may not have enough time to exactly 
reach a solution before a new input is needed by the system. 
In such cases, sub-optimal inputs to the control
system are used. 
These types of  
 \emph{feedback optimization} methods have been used 
 in various settings such as
 optimal power flow and human-in-the-loop control, among others~\cite{bernstein2019online,colombino2020online,dall2016optimal,lindemann2021learning,ospina2022feedback}. 


In this paper, we develop and analyze a distributed algorithm for multi-agent feedback optimization.
In centralized feedback optimization, 
all output measurements are fed into a single optimization algorithm, which computes
all inputs for the system. However, some modern control applications
consist of interacting decision-makers, such as buildings on the smart power grid,
and we therefore develop a
multi-agent feedback optimization framework. 
In this framework, different agents measure different system outputs and compute different system inputs, and
they communicate to collaborate. 
Many multi-agent systems face asynchrony in agents' communications, computations, and sensor measurements. For example, asynchronous communications can arise from 
adversarial jamming, 
asynchronous computations can stem from heterogeneous hardware, and 
asynchrony in sensor measurements can be due to intermittent feedback~\cite{sinopoli2004kalman}. 

Therefore, the goal of this article is to design a decentralized feedback optimization algorithm 
that enables a network of agents
to drive system outputs to track the 
time-varying solutions of time-varying optimization problems,
even when agents’ 
communications, computations, and sensor measurements are subject to asynchrony. 
In particular, we use a block-based gradient optimization algorithm. 
Asynchronous block-based algorithms were first established in seminal work in~\cite{tsitsiklis1986distributed,bertsekas1989convergence,bertsekas1989parallel}, and
recent developments have extended these results to constrained problems~\cite{hendrickson2022totally}, problems that satisfy the Polyak-Łojasiewicz condition~\cite{yazdani2021asynchronous}, and others~\cite{assran2020advances,ubl2021totally}. These successes motivate their use here as well. 
To the best of our knowledge, our work is the first to consider decentralized
feedback optimization with asynchrony
in computations, communications, and sensing.

To summarize, this paper makes the following contributions: 

\begin{itemize}
    \item We provide the first block-based asynchronous algorithm for distributed feedback optimization problems~{(Algorithm~\ref{alg:myAlg})}.
    \item We show that this algorithm converges toward the global minimizer for time-invariant and time-varying feedback optimization problems, and we derive a convergence rate~{(Theorems~\ref{thm:1} \& \ref{thm:2})}.
    \item We show that our algorithm 
    asymptotically tracks the global minimizer  
    to within an explicit error bound~{(Theorem~\ref{thm:3})}.
    \item We provide timing specifications for computations, communications, and sensor measurements to achieve desired network tracking performance {(Theorem~\ref{thm:4}, Corollary~\ref{cor:1})}.
    \item We empirically demonstrate the ability of our algorithm to track the solutions of feedback optimization problems by applying it to networks of 
    agents in two simulations~{(Section~\ref{sec:simulation})}.
\end{itemize}

Feedback optimization and time-varying optimization have been studied in the centralized setting
\cite{tang2018feedback,bernstein2019online}, as well as 
in the distributed setting~\cite{dall2016optimal,chang2019saddle,bolognani2014distributed}.
The most closely related works to the current article are~\cite{bernstein2018asynchronous,bastianello2023online,dall2016optimal,chang2019saddle,bolognani2014distributed}. Results in 
\cite{bernstein2018asynchronous} studied decentralized feedback optimization problems where agents' communications are subject to asynchrony, and 
\cite{bastianello2023online} considers time-varying optimization problems with asynchrony in agents' computations and communications, while
\cite{dall2016optimal} develops a primal-dual method utilizing feedback for the optimal power flow problem.
Efforts in 
\cite{chang2019saddle} developed a distributed saddle-flow algorithm to achieve consensus over a connected graph considering a feedback optimization problem, and finally
we mention that 
\cite{bolognani2014distributed} proposes a distributed feedback algorithm for the optimal power flow problem, and they consider asynchrony in agents' computations.
This article differs from all of these works in that we 
address asynchrony 
in agents' computations, communications, and sensor measurements simultaneously. 
This article extends our previous work~\cite{behrendt2021technical} on convex time-varying optimization problems,
and the current paper differs from~\cite{behrendt2021technical} because it considers the feedback optimization setting. 

The rest of this article is organized as follows.
Section~\ref{sec:problem} gives a problem statement, and
Section~\ref{sec:update} presents our asynchronous feedback optimization algorithm.
Convergence rates are derived in Section~\ref{sec:convergence}, 
Section~\ref{sec:simulation} presents simulations, 
and Section~\ref{sec:conclusion} concludes. \\

\noindent\textbf{Notation} Let~$\N=\{1,2,3,\dots \}$ denote the natural numbers and let~${\N_{0} = \N \cup \{ 0 \}}$ denote the non-negative integers. 
For~$N \in \N$, define~$[N] = \{1, \ldots, N\}$. 
We use~$\norm{\cdot}$ for the Euclidean norm. For $a \in \R^n$ and~$b \in \R^m$ we denote the stacking of these vectors as~$(a,b) \in \R^{n+m}$. 
The diameter of a compact set~$\X \subset \R^n$ is denoted~{$\diam ( \X ) \coloneqq  \sup_{x, y \in \mathcal{X}} \|x - y\|$}.
We use~$\Pi_{Z}$ for the Euclidean projection of a point
onto a closed, convex set~$Z$, i.e., $\Pi_Z[v] = \argmin_{z \in Z} \|v - z\|$.
We define~$\nabla_x \coloneqq \frac{\partial}{\partial x}$
and~$\nabla_y \coloneqq \frac{\partial}{\partial y}$. We denote~$I_n$ as the~$n \times n$ identity matrix
and~$\otimes$ is the Kronecker product. We also use~$\mathbb{1}_n$ as the~$n$-dimensional ones vector.

\section{Problem Formulation}\label{sec:problem}
This section states the problem
that is the focus of this paper. 

\subsection{Problem Statement}


Suppose the output of a dynamical system is given as
    ${y = h(x)}$,
where~$h:\R^n \rightarrow \R^m$ is a map from some controllable inputs~${x \in \X \subseteq \R^n}$ to outputs~$y \in \Y \subseteq \R^m$ for sets~$\X$ and~$\Y$. 
In this article, we consider systems where~$h$ is a linear map of the form 
    $y(k) = Cx(k)$, 
where~$C \in \R^{m \times n}$.
Our goal is to regulate the outputs of such a system 
to the solution of a time-varying optimization problem. 
We consider a network of agents doing so, 
and this setup models settings in which
computations of new inputs and/or measurements
of outputs are done in a parallelized fashion. 
That is, while the relationship~$y(k) = Cx(k)$ may hold,
measurements of each entry of~$y$ may happen onboard different embedded
sensors, and their embedded processors may compute new 
values for different entries of~$x$. 
This can be seen, e.g., in 
the real-time control of power distribution systems, which 
measures active and reactive powers from a network 
to regulate voltages and power 
flows~\cite{dall2016optimal}.
It can be difficult to synchronize agents' operations, and thus all agents are permitted to operate asynchronously as they compute 
new inputs, measure outputs, and communicate 
to work together. 
%

Formally, we consider problems of the following form. 
\begin{problem} \label{prob:main}
    Given~$f : \R^{n} \times \N_{0} \rightarrow \R$, $g : \R^{m} \times \N_{0} \rightarrow \R$, using~$N \in \N$ agents that asynchronously 
    compute, communicate, and measure outputs, drive~$x$ and~$y$ to 
    track the solution of
    \begin{align}
        \underset{x\in \X}{\textnormal{minimize}} & \ \J{x}{y}{\ell} \coloneqq \f{x}{\ell} + \g{y}{\ell} \\
        \textnormal{subject to} & \ y=Cx,
    \end{align}
    where~$t_{\ell} \in \T \coloneqq \{t_0, \ldots, t_T\}$
    and~$T \in \N$. 
    \hfill $\Diamond$
\end{problem}

Problem~\ref{prob:main} can model the scenario of sampling from a continuous-time objective function that 
was considered in~\cite{simonetto2016class,tang2022running}, though we model problems simply as occurring in discrete time. 

\begin{remark} \label{rem:disturbances}
Problem~\ref{prob:main} is an aggregated form of the problem 
that will be solved by~$N$ agents. However
no single agent will always know the most recent values
of all entries of~$y$ as it is written in Problem~\ref{prob:main}, and the same is true for~$x$. 
Instead, different agents will
measure different entries of~$y$, and different agents will compute 
new values for 
different entries of~$x$
(discussed in Section~\ref{sec:update}). Under asynchronous communications, agents will send and receive different
entries of~$x$ and~$y$ at different times, causing agents to have disagreeing local
copies of~$x$ and~$y$ onboard. 
These disagreements are disturbances in~$x$ and~$y$, in the sense that each
agent's local copy of~$x$ and~$y$ can be viewed as a perturbed version of the actual values of~$x$ and~$y$.
When~$x$ and~$y$ are used in local computations, their attendant perturbations enter these computations as well. 
Feedback optimization has been shown to be robust to various forms of perturbations, and
in Section~\ref{sec:convergence}, we show that it is robust to the perturbations that result
from asynchrony as well. \hfill $\blacklozenge$
\end{remark}

\subsection{Assumptions on Problem~\ref{prob:main}}
We make the following assumptions about Problem \ref{prob:main}.

\begin{assumption} \label{ass:convex}
For all~$\ind{\ell}  \in \T$, the functions~$\f{\cdot}{\ell}$ and~$\g{\cdot}{\ell}$ are twice continuously differentiable and~{$p$-strongly} convex. \hfill $\lozenge$ 
\end{assumption}

Assumption~\ref{ass:convex} guarantees the existence and continuity of both the gradient and Hessian of~$\J{\cdot}{\cdot}{\ell}$. Additionally, it implies that~$\J{\cdot}{\cdot}{\ell}$ is~$p$-strongly convex for all~$\ind{\ell} \in \T$. 
\begin{assumption} \label{ass:setConstraint}
    The constraint set~$\X$ can be decomposed via
        ${\X = \X_{1} \times \X_{2} \times \dots \times \X_{N}}$,
    where~$\X_{i} \subseteq \R^{n_i}$,~$n_i \in \N$, is non-empty, compact, and polyhedral for all~$i \in [N]$. \hfill $\lozenge$ 
\end{assumption}

Assumption~\ref{ass:setConstraint} permits constraints such as box constraints, which
are common in multi-agent problems. 
Because~$y = Cx$ we also have~$y \in \Y$, where~$\Y = \{y = Cx : x \in \X\}$,
and under Assumption~\ref{ass:setConstraint} the set~$\Y$
is non-empty, compact, and convex. 

We decompose~$x\in\X$ via~$x = \begin{bmatrix} x^T_1 \hdots x^T_N  \end{bmatrix}^T$, where for all ${i \in [N]}$ we have~$x_{i} \in \X_{i} \subseteq \R^{n_i}$ and~$n= \sum^N_{i=1} n_i$. We decompose~$y \in \Y$ via~$y = \begin{bmatrix} y^T_1 \hdots y^T_N  \end{bmatrix}^T$, where~$y_{i} \in \R^{m_i}$ and~$m= \sum^N_{i=1} m_i$. For notational simplicity,
we consider the case where every agent measures at least one output,
i.e., at least one entry of~$y$, though all of our developments directly apply to problems in which this is not the case.
Under this decomposition of~$y$, agent~$i$ measures the block of outputs~$y_i$
and performs computations to determine subsequent values of~$x_i$. 
We note that, for all~$i \in [N]$, 
Assumption~\ref{ass:setConstraint} allows agent~$i$ to project values of~$x_i$ onto~$\X_i$, and,
when done by all agents, this projection ensures that~$x\in\X$. 

It will be helpful in the forthcoming analysis to partition the matrix~$C$ as
    $C = \begin{bmatrix}
        C_1 & C_2 & \dots & C_N
    \end{bmatrix}$,
where we have defined the matrix~$\Ci \in \R^{m \times n_i}$ for all~$i \in [N]$. 
Here,~$C_1$ is the first~$n_1$ columns of~$C$, then~$C_2$
is the next~$n_2$ columns, etc. 
We also denote the rows of~$C$ via
    $C = \begin{bmatrix}
        C^T_{1*} & C^T_{2*} & \dots & C^T_{N*}
    \end{bmatrix}^T$,
where~$C_{1*}$ is the first~$m_1$ rows, then~$C_{2*}$
is the next~$m_2$ rows, etc.

 Assumptions~\ref{ass:convex} and~\ref{ass:setConstraint} ensure the existence and uniqueness of 
 the minimizer of~$J(\cdot, \cdot; t_{\ell})$.
 For~$t_{\ell} \in \T$ we denote this minimizer as 
    \begin{equation} \label{eq:solution}
        \x{*}{}{\ind{\ell}} \coloneqq \underset{x \in \X}{\argmin} \ \f{x}{\ell} + \g{Cx}{\ell}, \quad 
        \y{*}{}{\ind{\ell}} = C \x{*}{}{\ind{\ell}}.
    \end{equation}
    Assumptions~\ref{ass:convex} and~\ref{ass:setConstraint} also give the following lemma.

\begin{lemma} [Error Bound Condition] \label{lem:errorBound}
Let Assumptions~\ref{ass:convex} and~\ref{ass:setConstraint} hold. Then
    for every~$\varpi > 0$ and for each~$\ind{\ell} \in \T$, there exist~$\upsilon, \lambda > 0$ such that for all~$x \in \X$ and~$y \in \Y$ with (i)~$\J{x}{y}{\ell} \leq \varpi$ and (ii)~$\big\Vert x - \Pi_\X \big[ \x{}{}{} - \gamma_\ell \nabla_x \J{\x{}{}{}}{\y{}{}{}}{\ell} \big] \big\Vert \leq \upsilon$, we have the upper bounds 
        $\norm{x - \x{*}{}{\ind{\ell}}} \leq \lambda \big\Vert x - \Pi_\X \big[ \x{}{}{} - \nabla_x \J{\x{}{}{}}{\y{}{}{}}{\ell} \big] \big\Vert$
    and
    \begin{equation}
        \norm{x - \x{*}{}{\ind{\ell}}} \leq \lambda \max\{ 1,\gamma^{-1}_\ell \} \big\Vert x - \Pi_\X \big[ \x{}{}{} -  \gamma_\ell \nabla_x \J{\x{}{}{}}{\y{}{}{}}{\ell} \big] \big\Vert \label{eq:15}
    \end{equation}
    for all~$\gamma_\ell > 0$. \hfill $\blacksquare$
\end{lemma}

Lemma~\ref{lem:errorBound} (in its original form without outputs) holds for a number of problem classes as discussed in~\cite[Section 2]{tseng1991rate,pang1987aposteriori,robinson1981some,luo1992linear,drusvyatskiy2018error,drusvyatskiy2018error,zhang2017restricted}. 
In this article,~$\J{\cdot}{\cdot}{\ind{\ell}}$ satisfies the error bound condition for each~$\ind{\ell} \in \T$ because it is strongly convex (by Assumption~\ref{ass:convex}) and defined over a polyhedral set (by Assumption~\ref{ass:setConstraint}),
which is established in~\cite{pang1987aposteriori}. 
The convergence proofs of this article 
therefore use Lemma~\ref{lem:errorBound} and 
draw in part on the work in~\cite{tseng1991rate}, which derives a linear convergence rate for asynchronous projected gradient iterations for a class of time-invariant problems (that are not feedback optimization)
that satisfy Lemma~\ref{lem:errorBound}.

For the time evolution of Problem~\ref{prob:main}, we assume the following. 
 \begin{assumption} \label{ass:time2}
     For all~$x \in \X$,~$y \in \Y$, and~$\ind{\ell} \in \T$, there exist~$L_t,\Delta , \sigma_{\ell} > 0$ such that 
     (i)   ${\big\Vert \x{*}{}{\ind{\ell+1}} - \x{*}{}{\ind{\ell}} \big\Vert \leq \sigma_{\ell+1}}$, 
     (ii)~$\big\vert \J{x}{y}{\ell + 1} - \J{x}{y}{\ell} \big\vert \leq L_t \big\vert \ind{\ell + 1} - \ind{\ell} \big\vert$, and 
     (iii)~${\big\vert \ind{\ell+1} - \ind{\ell} \big\vert \leq \Delta}$.    
     \hfill $\lozenge$
 \end{assumption}
 
The first condition in Assumption~\ref{ass:time2} 
ensures that successive minimizers 
are not arbitrarily far apart. 
Without such an assumption it may be impossible for an algorithm to track 
time-varying solutions, and therefore we enforce this condition here
to ensure that Problem~\ref{prob:main} is solvable. 
The second condition 
establishes a Lipschitz continuous-like condition in time.
The last condition in Assumption~\ref{ass:time2} simply states that the time between changes in the objective function is bounded.

\section{Asynchronous Update Law} \label{sec:update}
In this section we develop the proposed asynchronous algorithm we use to 
drive agents' outputs to
track the solution of Problem~\ref{prob:main} over time. 
In this section, we use the term ``operations" to collectively refer to agents' computations, communications, and output measurements.

\subsection{Timescale Separation}
The objective in Problem~\ref{prob:main} is indexed by the discrete time index~$\ind{\ell}$, and
we index agents' operations over a different discrete time index, namely~$k$, because 
the timing of agents' operations can differ from the timing of changes
in their objective functions. In fact,~\cite{simonetto2017decentralized}
notes that for correction-only algorithms, namely, algorithms that do not predict future objective functions,
some timescale separation is required between the changes in objective functions and the agents'
operations. 
That is, between the change from~$\ind{\ell}$ to~$\ind{\ell +1}$ there must be some non-zero number
of ticks of~$k$. We consider a correction-only algorithm since it may be difficult to predict
discrete jumps in agents' objective function, and
we therefore assume the following.

\begin{assumption} \label{as:timescale}
In Problem~\ref{prob:main}, for each~$\ind{\ell} \in \T$, there are~$\kappa_{\ell} \geq 1$ ticks of~$k$ when minimizing~$\J{\cdot}{\cdot}{\ell}$. 
\hfill $\lozenge$
\end{assumption}

This assumption is one of technical feasibility rather than conveience; 
without Assumption~\ref{as:timescale}, agents may not be able to track solutions
at all, and this assumption at least makes it possible to track a solution
with bounded error. Such tracking is not guaranteed, and we still must
devise an algorithm and characterize its performance.

We use~$\eta_\ell \in \N$ to denote the total number of ticks of~$k$ that have elapsed from~$t_0$ to the moment  before~$\ind{\ell}$ increments to~$\ind{\ell + 1}$, i.e.,
 \begin{equation} \label{eq:etadef}
 \eta_\ell = \sum^{\ell}_{i=0} \kappa_i.
\end{equation}

\subsection{Formal Algorithm Statement}
We consider a block-based gradient projection algorithm with asynchronous computations, communications, and output measurements to make~$x$ and~$y$ approximately 
track~$\big\{\big(x^*(t_{\ell}), y^*(t_{\ell})\big)\big\}_{t_{\ell} \in \T}$. 
Each agent updates only a subset of the entries of
the input vector,~$x$, and 
each agent measures only a subset of the entries of the output vector,~$y$. 
Over time, each agent locally computes new values for its entries of~$x$ and then
communicates these new values and values of~$y$ that it has measured
to other agents. 

Asynchrony implies that agents receive different information at different times, and thus we expect them to have differing values for network inputs and outputs onboard. 
At any time~$k$, agent~$i$ has a local copy of the network input and output vectors, denoted as~$\x{i}{}{k}$ and~$\y{i}{}{k}$, respectively. 
Due to asynchrony, we allow~$\x{i}{}{k} \neq \x{j}{}{k}$ and~$\y{i}{}{k} \neq \y{j}{}{k}$ for~$j \neq i$. Within the vector~$\x{i}{}{k}$, agent~$i$ computes new values only for its own sub-vector of inputs, which is~$\x{i}{i}{k} \in \R^{n_i}$. 
Similarly, within the vector~$\y{i}{}{k}$, agent~$i$ measures only a sub-vector of outputs, denoted~$\y{i}{i}{k} \in \R^{m_i}$. 
At any time~$k$, agent~$i$ has onboard (possibly old) values for agent~$j$'s sub-vector of the inputs and outputs, denoted by~$\x{i}{j}{k} \in \R^{n_j}$ and~$\y{i}{j}{k} \in \R^{m_j}$, respectively. 
These values onboard agent~$i$ only change when agent~$i$ receives
a  communication from agent~$j$. In particular, agent~$i$ does
not perform any computations on~$\x{i}{j}{k}$ and does
not measure~$\y{i}{j}{k}$ at any point in time; only agent~$j$
does these operations. 

At time~$k$, if agent~$i$ computes an update to~$\x{i}{i}{k}$,
then it performs these computations with its onboard
input vector~$\x{i}{}{k}$ and onboard output vector~$\y{i}{}{k}$ 
because these are all that it has access to. 
As noted in the preceding paragraph, the entries of~$\x{i}{j}{k}$ and~$\y{i}{j}{k}$ for~$j \! \neq \! i$ are obtained by communications from agent~$j$, which may be subject to delays. Therefore, agent~$i$ may (and often will) compute updates to~$\x{i}{i}{k}$ using outdated information from other agents. 

To formalize an algorithm statement, 
let~$\K^i$ be the set of times at which agent~$i$ computes an update to~$\x{i}{i}{}$. Similarly, let~$\M^i$ be the set of times at which agent~$i$ takes measurements of~$\y{i}{i}{}$. 
Let~$\C^j_i$ be the set of times at which agent~$i$ receives transmission of~$\x{i}{j}{}$ and~$\y{i}{j}{}$ from agent~$j$; due to communication delays, these transmissions can be received at some time after they are sent, and they can be received at different times by different agents. 
We emphasize that the sets~$\K^i, \ \M^i$, and~$\C^j_i$ are only defined to simplify discussion; agents do not know (and do not need to know)~$\K^i, \ \M^i$, and~$\C^j_i$. 

We define~$\tau^i_j(k)$ to be the time at which agent~$j$ originally computed the value of~$\x{i}{j}{k}$ that agent~$i$ has onboard at time~$k$. 
We define~$\mu^i_j(k)$ to be the time at which agent~$j$ originally measured the value of~$\y{i}{j}{k}$ that agent~$i$ has onboard at time~$k$. 
Using this notation, at any time~$k$, agent~$i$ stores onboard
\begin{align}
    \x{i}{}{k} &= \left( \x{1}{1}{\tau^i_1(k)}^T,\dots,\x{i}{i}{k}^T,\dots, \x{N}{N}{\tau^i_N(k)}^T \right)^T  \label{eq:3} \\
    \y{i}{}{k} &= \left( \y{1}{1}{\mu^i_1(k)}^T,\dots,\y{i}{i}{k}^T,\dots, \y{N}{N}{\mu^i_N(k)}^T \right)^T. \label{eq:4}
\end{align}
In Section~\ref{sec:convergence}, we will also analyze the ``true'' state of the network, 
\begin{align}
    \x{}{}{k} &= \left( \x{1}{1}{k}^T,\x{2}{2}{k}^T,\dots,\x{N}{N}{k}^T \right)^T \label{eq:85} \\
    \y{}{}{k} &= \left( \y{}{1}{k}^T,\y{}{2}{k}^T,\dots,\y{}{N}{k}^T \right)^T, \label{eq:86}
\end{align} 
where~$\y{}{i}{k} = C_{i*} \x{}{}{k}$. Here~$x(k)$ and~$y(k)$ contain all of the
most recent values of inputs and outputs; no agent knows these vectors and they
are used only for analysis. 

\begin{remark}
As noted in Remark~\ref{rem:disturbances}, asynchrony causes disturbances in agent~$i$'s local copy of~$x$ and~$y$.
In particular, the disturbances in~$x^i(k)$ and~$y^i(k)$ are~$x^i(k) - x(k)$ and~$y^i(k) - y(k)$, respectively.
These disturbances are not known to agent~$i$, though
we show in Section~\ref{sec:convergence} that feedback optimization is robust
to them, despite not having an explicit model for them, and this robustness
is in line with the existing feedback optimization literature~\cite{colombino2019towards,bernstein2019online}. 
\hfill $\blacklozenge$
\end{remark}

We assume that computation, communication, and measurement delays are bounded, which is called \emph{partial asynchrony}. 

\begin{assumption}[Partial Asynchrony] \label{ass:partialAsynch}
Let~$\K^i$ be the set of times at which agent~$i$ computes an update to~$\x{i}{i}{}$, 
and let~$\M^i$ be the set of times at which agent~$i$ takes measurements of~$\y{i}{i}{}$. 
Let~$\tau^i_j(k)$ be the time at which agent~$j$ originally computed the value of~$\x{i}{j}{k}$ that agent~$i$ has onboard at time~$k$,
and let~$\mu^i_j(k)$ to be the time at which agent~$j$ originally measured the value of~$\y{i}{j}{k}$ that agent~$i$ has onboard at time~$k$. 
Then there exists~$B \in \N$ such that: 
\begin{enumerate}
\item For every $i \in [N]$ and for every $k\geq0$, at least one of the elements
of the set $\left\{ k,k+1,\dots,k+B-1\right\} $ belongs to $\K^{i}$\label{as:pa1}
\item For every $i \in [N]$ and for every $k\geq0$, at least one of the elements
of the set $\left\{ k,k+1,\dots,k+B-1\right\} $ belongs to $\M^{i}$\label{as:pa2}
\item It holds that
$\max\left\{ 0,k-B+1\right\} \leq\tau_{j}^{i}\left(k\right)\leq k$
 for all $i \in [N]$, $j \in [N]$, and $k\geq0$.\label{as:pa3}
\item It holds that
$\max\left\{ 0,k-B+1\right\} \leq\mu_{j}^{i}\left(k\right)\leq k$
 for all ${i \in [N]}$, $j \in [N]$, and $k\geq0$.\label{as:pa4} \hfill $\lozenge$ 
\end{enumerate}
\end{assumption}

Assumptions~\ref{ass:partialAsynch}.1 and \ref{ass:partialAsynch}.2 ensure each agent updates their assigned decision variables and measures their assigned outputs at least once every~$B$ time steps. Assumptions~\ref{ass:partialAsynch}.3 and \ref{ass:partialAsynch}.4 ensure that agents communicate those updates and measurements at least once every~$B$ time steps. 
We will use this assumption to prove that for each~$\ind{\ell} \in \T$ our algorithm will 
provably make progress toward~$\big( \x{*}{}{\ind{\ell}}, \y{*}{}{\ind{\ell}} \big)$ across each
interval of~$B$ time steps. 
For simplicity, for each~$\ind{\ell} \in \T$
we consider~$\kappa_\ell = r_\ell B$ for some~$r_\ell \in \N_0$. 
Furthermore, we adopt the convention that~$\kappa_{-1}=0$ and~$\eta_{-1}=0$.

We seek to develop a distributed 
feedback optimization 
algorithm that is robust to asynchrony.
It has been shown that gradient-based methods are robust
to asynchrony for static optimization problems~\cite{tsitsiklis1986distributed,bertsekas1989convergence,bertsekas1989parallel}, and we 
therefore we propose the update law 
\begin{align}
\x{i}{i}{k+1}= & \begin{cases} 
\zeta_i\big(x^i(k), y^i(k)\big) & k\in \K^{i}\\
x_{i}^{i}(k) & k\notin \K^{i}
\end{cases}\label{eq:update}\\
\x{i}{j}{k+1}= & \begin{cases}
\x{j}{j}{\tau_{j}^{i}(k)} & k\in \C_{i}^{j}\\
\x{i}{j}{k} & k\notin \C_{i}^{j},
\end{cases}
\end{align}
where
\begin{multline} \label{eq:zetai}
\zeta_i\big(x^i(k), y^i(k)\big) \\ = \Pi_{\X_{i}}\Big[\x{i}{i}{k} - \gamma_{\ell}\Big(\nabla_{x_{i}} \f{\x{i}{}{k}}{\ell} + \Ci^{T}\nabla_{y} \g{\y{i}{}{k}}{\ell} \Big)\Big]
\end{multline}
and~$\gamma_{\ell} > 0$ is the step size used to minimize~$J(\cdot, \cdot; t_{\ell})$. 
In addition to agents' computations, we consider the output law
\begin{align}
\y{i}{i}{k+1}= & \begin{cases}
\y{}{i}{k} & k\in \M^{i}\\
\y{i}{i}{k} & k\notin \M^{i}
\end{cases}\\
\y{i}{j}{k+1}= & \begin{cases}
\y{j}{j}{\mu_{j}^{i}(k)} & k\in \C_{i}^{j}\\
\y{i}{j}{k} & k\notin \C_{i}^{j},
\end{cases}
\end{align}
where we define~$y_i(k) = C_{i*}x(k)$. 

This update law only requires agents to perform computations with
their local copies of the network inputs and outputs. 
Algorithm~\ref{alg:myAlg} provides pseudocode for agents' update law, and
the remainder of this paper focuses on analyzing the execution of Algorithm~\ref{alg:myAlg}.

\begin{algorithm} 
    \caption{Asynchronous Feedback Optimization Algorithm}
    \label{alg:myAlg}
	\SetAlgoLined
	\KwIn{$\x{i}{j}{0}$ and~$\y{i}{j}{0}$ for all~$i,j \in [N]$}
	\For{$\ind{\ell} \in \T$}
	{
		\For{$k=\eta_{\ell-1}+1:\eta_{\ell}$}
		{
			\For{i=1:N}
			{
			    \For{j=1:N}
			    {
    				\If{$k \in C^{j}_i$}{    			
                        ~\phantom{} $\x{i}{j}{k+1} = \x{j}{j}{\tau^i_j(k)}$ \\
                        $\y{i}{j}{k+1} = \y{j}{j}{\mu^i_j(k)}$
    				}
    				\Else{
    				    ~\phantom{} $\x{i}{j}{k+1} = \x{i}{j}{k}$ \\
                        $\y{i}{j}{k+1} = \y{i}{j}{k}$
    				}
    			}
				\If{k $\in \K^{i}$ }{
					~$\x{i}{i}{k+1} = \zeta_i\big(x^i(k), y^i(k)\big)$
				}
				\Else{
				~$\x{i}{i}{k+1} = \x{i}{i}{k}$ 
				}
                \If{k $\in \M^{i}$ }{
					~$\y{i}{i}{k+1} = \y{}{i}{k}$
				}
				\Else{
				~$\y{i}{i}{k+1} = \y{i}{i}{k}$ 
				}
			}			
		}	
	}
\end{algorithm}

\section{Convergence of Asynchronous Feedback Optimization Algorithm}
\label{sec:convergence}
In this section we prove the approximate convergence of Algorithm~\ref{alg:myAlg}. 
First, we establish the following properties of the objective function.
For all~$\ind{\ell} \in \T$, from Assumptions \ref{ass:convex}-\ref{ass:setConstraint} both~$\nabla_x^2 \f{\cdot}{\ell}$ and~$\nabla_y^2 \g{\cdot}{\ell}$ are continuous and
both~$\X$ and~$\Y$ are compact. 
Therefore~$\nabla_x \f{\cdot}{\ell}$ and~$\nabla_y \g{\cdot}{\ell}$ are Lipschitz on~$\X$ and~$\Y$, respectively. For each~$\ind{\ell} \in \T$, we use~$\Lx$ and~$\Ly$ to denote upper bounds on~$\norm{\nabla_x^2 \f{\cdot}{\ell}}$ and~$\norm{\nabla_y^2 \g{\cdot}{\ell}}$ over~$\X$ and~$\Y$, respectively. These~$\Lx$ and~$\Ly$ are the Lipschitz constants of~$\nabla_x \f{\cdot}{\ell}$ and~$\nabla_y \g{\cdot}{\ell}$ over~$\X$ and $\Y$, respectively. Thus, for each~${\ind{\ell} \in \T}$, and for all~$x_1,x_2 \in \X$ and~$y_1,y_2 \in \Y$, we have 
\begin{align}
    \big\Vert \nabla_x \f{x_1}{\ell} - \nabla_x \f{x_2}{\ell} \big\Vert &\leq \Lx \big\Vert x_1 - x_2 \big\Vert\\
    \big\Vert \nabla_y \g{y_1}{\ell} - \nabla_y \g{y_2}{\ell} \big\Vert &\leq \Ly \big\Vert y_1 - y_2 \big\Vert.
\end{align}
By an analogous argument, for all~$x \in \X$,~$y \in \Y$, and~$\ind{\ell} \in \T$, there exist~$M_{x,\ell}>0$ and~$M_{y,\ell}>0$ such that
\begin{equation} \label{eq:gradbound}
    \Vert \nabla_x \f{x}{\ell} \Vert \leq M_{x,\ell} \qquad \textnormal{ and } \qquad
    \Vert \nabla_y \g{y}{\ell} \Vert \leq M_{y,\ell}.
\end{equation}
Since~$\nabla \J{\cdot}{\cdot}{\ell}$ is continuous over the compact set
$\X \times \Y$, it is bounded on~$\X \times \Y$. Thus, 
the function~$\J{\cdot}{\cdot}{\ell}$ is Lipschitz for each~$\ind{\ell} \in \T$,
and there exists a Lipschitz constant~${\Lj > 0}$ such that
    $\big\vert\J{x_1}{y_1}{\ell} - \J{x_2}{y_2}{\ell} \big\vert \leq \Lj \norm{ \big(x_1,y_1 \big) - \big(x_2,y_2 \big)}$
for all~$x_1, x_2 \in \X$ and~$y_1, y_2 \in \Y$. 
By a similar argument, continuity of~$\nabla_x^2 \J{\cdot}{\cdot}{\ell}$ 
over~$\X \times \Y$ implies that~$\nabla_x \J{\cdot}{\cdot}{\ell}$
is Lipchitz, and
for each~$\ind{\ell} \in \T$ 
there exists a Lipschitz constant~$\Lt > 0$ such that
    $\norm{\nabla_x \J{x_1}{y_1}{\ell} - \nabla_x \J{x_2}{y_2}{\ell}}\leq \Lt \norm{ \big(x_1,y_1 \big) - \big(x_2,y_2 \big)}$
for all~$x_1,x_2 \in \X$ and~$y_1,y_2 \in \Y$.  

To keep track of computations, we define 
\begin{align}
    \s{i}{k}{\ell} & \coloneqq \begin{cases}
                    \x{i}{i}{k+1} - \x{i}{i}{k} & k\in \K^{i}\\
                    0 & k \notin \K^{i}, \label{eq:81}
\end{cases}
\end{align} 
and to track output measurements we define~$\q{i}{k}$ as
\begin{align}
    \q{i}{k} & \coloneqq  \begin{cases}
                    \y{i}{i}{k+1} - \y{i}{i}{k} & k\in \M^{i}\\
                    0 & k \notin \M^{i},
\end{cases} 
\end{align}
where we note~$\y{i}{i}{k+1} - \y{i}{i}{k}=\y{}{i}{k} - \y{}{i}{\mu_i^i(k)}$
for~$k \in \M^{i}$.

We concatenate terms in~$\s{}{k}{\ell} \coloneqq \big( \s{1}{k}{\ell}^T, \dots, \s{N}{k}{\ell}^T \big)^T \in \R^n$, and ~$\q{}{k} \coloneqq \big( \q{1}{k}^T, \dots, \q{N}{k}^T \big)^T \in \R^m$.
We also define the following terms for all~$\ind{\ell} \in \T$ and for all~$k$ such that~$\eta_{\ell -1} \leq k \leq \eta_\ell$:
\begin{align}
    \at{k}{\ell} &\coloneqq \J{\x{}{}{k}}{\y{}{}{k}}{\ell} - \J{\x{*}{}{\ind{\ell}}}{\y{*}{}{\ind{\ell}}}{\ell} \\
    \bt{k}{\ell} &\coloneqq \sum_{\tau=k-B }^{k-1} \norm{\s{}{\tau}{1}}^2, \qquad \dt{k}{\ell} \coloneqq \sum_{\tau=k-B }^{k-1} \norm{\q{}{\tau}{}}^2. \label{eq:104}
\end{align}

The next theorem establishes a convergence rate for 
the first, static objective function that agents minimize.
After that, we will 
extend our analysis to a sequence of time-varying objectives 
in Theorem~\ref{thm:2}. 

\begin{theorem} \label{thm:1}
    Let Assumptions~\ref{ass:convex}, \ref{ass:setConstraint}, and \ref{ass:partialAsynch} hold. 
    Then, for fixed~$t_0$, 
    a stepsize~$\gamma_0 \in (0,\gamma_{\max,0})$,
    iterates~$\x{}{}{k}$ and~$\y{}{}{k}$ as defined in~\eqref{eq:85} and~\eqref{eq:86}, and any~$r_0 \in \mathbb{N}$, the sequence~$\{ \x{}{}{k}, \y{}{}{k} \}_{k \in \N_0}$ generated by~$N$ agents executing Algorithm~\ref{alg:myAlg} satisfies
    \begin{align}
        \at{r_0 B}{0} \leq a_0 \pt{r_0 - 1}{0} \label{eq:t1_alpha} \\
        \bt{r_0 B}  &\leq b_0 \pt{r_0 - 1}{0} \label{eq:t1_beta} \\
        \dt{r_0 B} &\leq d_0 \pt{r_0-1}{0}, \label{eq:t1_delta}
    \end{align}
where~$D_0$ and~$E_0$ are from~$\eqref{eq:bigone_dande}$,~$F_0$ is from~\eqref{eq:bigone_f},~$G_0$ is from~\eqref{eq:bigone_g},
$\gamma_{\max, 0}$ is from~\eqref{eq:bigone_gmax},~$c_0$,~$\rho_0$,~$b_0$, and~$d_0$ are from~\eqref{eq:bigone_candrhoandbandd}, 
    ${a_0 \!=\! \max \!\Bigg\{\! L_{J,0} \big( 1 \!+\! \norm{C} \big) \diam(\X), 
    \!\frac{8 E_0 \Big( \frac{G_0}{F_0} + \frac{E_0}{D_0} \Big) F_0}{D_0} B \diam\left(\X\right)^{2} \!\!\Bigg\}}$
    and $b_0 = \frac{D_0}{8 E_0 \Big( \frac{G_0}{F_0} + \frac{E_0}{D_0} \Big) F_0} a_0$. 
\end{theorem}
\emph{Proof:} See Appendix~B. \hfill $\blacksquare$

\begin{figure*}[t!]
\begin{align}
    D_\ell &=  \frac{  2 - \gamma_\ell \big(\big( 1+ B \big) \Lx + \big( 1 + B N  \big) \norm{C}^2 \Ly \big) }{2}, \qquad\qquad  E_\ell =   N B \frac{ \Lx + \Ly N \norm{C}^2}{2} \label{eq:bigone_dande} \\
    F_\ell &= \frac{1}{2} 
    \bigg((1 + \lambda^2) \Big[ 36 B^3 \Vert C \Vert^6 L_\ell^2 L_{y,\ell}^2 N^2 m 
    + 72 B^3 \Vert C \Vert^4 L_\ell^2 L_{x,\ell} L_{y,\ell} N^2 m 
    + 36 B L_\ell^2 L_{x,\ell}^2 N^2 
    + 36 B^3 \Vert C \Vert^2 L_\ell^2 L_{x,\ell}^2 N^2 m \\
    &+ 36 B \Vert C \Vert^4 L_\ell^2 L_{y,\ell}^2 N^2 
    + 18 \Vert C \Vert^2 L_{x,\ell} L_{y,\ell} N 
    + 72 B \Vert C \Vert^2 L_\ell^2 L_{x,\ell} L_{y,\ell} N^2 
    + 9 \Vert C \Vert^4 L_{y,\ell}^2 N\Big] 
    + 3 B^2 \Vert C \Vert^6 L_\ell^2 L_{y,\ell}^2 N^2 m \\
    &+ \Vert C \Vert^4 
    \big[72 B^3  L_\ell^2 L_{y,\ell} N^2 m 
    + 6 B^2 L_\ell^2 L_{x,\ell} L_{y,\ell} N^2 m 
    + 6 B^2 L_\ell^2 L_{y,\ell} N^2 m 
    + 3 L_\ell^2 L_{y,\ell}^2 N^2 
    + 3 B^2 L_{y,\ell}^2 N m\big] \\
    &+ \Vert C \Vert^2 \Big[ 96 B^3  L_\ell^2 N^2 m
    + 72 B^3  L_\ell^2 L_{x,\ell} N^2 m 
    + 72 B  L_\ell^2 L_{y,\ell} N^2
    + 60 B^3  L_\ell^2 N^2  \lambda^2 m
    + 18  L_{y,\ell} N 
    + 8 B^2  L_\ell^2 N^2 m  \\
    &+ 6 B^2  L_\ell^2 L_{x,\ell} N^2 m 
    + 6  L_\ell^2 L_{x,\ell} L_{y,\ell} N^2  
    + 6  L_\ell^2 L_{y,\ell} N^2
    + 3 B^2  L_\ell^2 L_{x,\ell}^2 N^2 m \Big]
    + 3 L_\ell^2 L_{x,\ell}^2 N^2 
    + 96 B L_\ell^2 N^2  
    + 6 L_\ell^2 L_{x,\ell} N^2  \\
    &+ 8 L_\ell^2 N^2 
    + 60 B L_\ell^2 N^2  \lambda^2 
    + 72 B L_\ell^2 L_{x,\ell} N^2
    + 12 L_{x,\ell}^2 N 
    + 9 L_{x,\ell}^2 N  \lambda^2
    + 18 L_{x,\ell} N 
    + 15 N  \lambda^2 
    + 24 N 
    + 2 \bigg) \label{eq:bigone_f} \\
    G_\ell &= \frac{N}{2} 
    \bigg( (1 + \lambda^2) \Big[
    72 B^3 \Vert C \Vert^4 L_\ell^2 L_{x,\ell} L_{y,\ell} N m 
    + 72 B \Vert C \Vert^2 L_\ell^2 L_{x,\ell} L_{y,\ell} N  
    + 36 B^3 \Vert C \Vert^6 L_\ell^2 L_{y,\ell}^2 N m 
    + 36 B^3 \Vert C \Vert^2 L_\ell^2 L_{x,\ell}^2 N m \\
    &+ 36 B \Vert C \Vert^4 L_\ell^2 L_{y,\ell}^2 N 
    + 36 B L_\ell^2 L_{x,\ell}^2 N
    \Big]
    + 3 B^2 \Vert C \Vert^6 L_\ell^2 L_{y,\ell}^2 N m
    + \Vert C \Vert^4 \Big[ 
    72 B^3  L_\ell^2 L_{y,\ell} N m 
    + 6 B^2  L_\ell^2 L_{x,\ell} L_{y,\ell} N m \\
    &
    + 6 B^2  L_\ell^2 L_{y,\ell} N m 
    + 3 B^2  L_{y,\ell}^2 m 
    + 3  L_\ell^2 L_{y,\ell}^2 N \Big] 
    + \Vert C \Vert^2 \Big[
    96 B^3  L_\ell^2 N m 
    + 72 B^3  L_\ell^2 L_{x,\ell} N m 
    + 72 B  L_\ell^2 L_{y,\ell} N
    + 60 B^3  L_\ell^2 N  \lambda^2 m \\ 
    &+ 8 B^2  L_\ell^2 N m 
    + 6 B^2  L_\ell^2 L_{x,\ell} N m 
    + 6  L_\ell^2 L_{x,\ell} L_{y,\ell} N 
    + 6  L_\ell^2 L_{y,\ell} N 
    + 3 B^2  L_\ell^2 L_{x,\ell}^2 N m
    + B  L_{y,\ell} N  \Big]
    + 96 B L_\ell^2 N 
    + 72 B L_\ell^2 L_{x,\ell} N \\
    &+ 60 B L_\ell^2 N  \lambda^2 
    + 8 L_\ell^2 N 
    + 6 L_\ell^2 L_{x,\ell} N 
    + 3 L_{x,\ell}^2 
    + 3 L_\ell^2 L_{x,\ell}^2 N 
    + B L_{x,\ell} 
    \bigg) \label{eq:bigone_g} \\
    a_\ell &= a_{\ell-1} \pt{r_{\ell-1} -1}{\ell-1} \!\!+ 2 \Delta L_t \!+\! L_{J,\ell} \sigma_{\ell} \big( 1 \!+\! \big\Vert C \big\Vert \big) \big(  M_{x,\ell} \!+\!  M_{y,\ell} \Vert C \Vert \big)  B \diam\left( \X \right) + \frac{8 E_\ell\left(\frac{G_\ell}{F_\ell} \!+\! \frac{E_\ell}{D_\ell}\right) F_\ell B^2 \diam \left(\X\right)^{2} \big( L_{x,\ell} \!+\! L_{y,\ell} \Vert C \Vert^2  \big)}{2 D_\ell} \label{eq:bigone_a} \\
    \gamma_{\max,\ell} &= \min \bigg\{\frac{2}{\big( 3N+1 \big) B \Lx + \big( 3N^2+1 \big) B \norm{C}^2 \Ly}, \frac{2}{\big( 1+ B \big) \Lx + \big( 1 +B N \big) \norm{C}^2 \Ly}, 
\frac{D_{\ell}}{E_{\ell}},
\bigg(\frac{G_{\ell}}{F_{\ell}} + \frac{E_{\ell}}{D_{\ell}} \bigg)^{-1},
\frac{1}{2c_\ell}, \\
&\qquad\qquad\qquad\qquad\qquad\qquad\qquad\qquad\frac{D_{\ell}}{8 F_{\ell} \Big( \frac{G_{\ell}}{F_{\ell}} + \frac{E_{\ell}}{D_{\ell}} \Big) c_\ell},
\frac{ \frac{a_\ell}{b_\ell} + 2E_\ell + D_{\ell} c_\ell - \sqrt{\big( \frac{a_\ell}{b_\ell} + 2E_\ell + D_{\ell} c_\ell \big)^2 - 4 D_{\ell} E_{\ell} c_\ell}}{2 E_{\ell} c_\ell}, 
\frac{1}{2} \bigg\}, \label{eq:bigone_gmax} \\
c_\ell &= \frac{D_\ell}{2 F_\ell + 2 D_\ell} \in \left(0,\frac{1}{2} \right), \qquad \pt{}{\ell} = 1 - \gamma_\ell c_\ell \in (0,1), \qquad b_\ell = B \diam \left(\X\right)^{2}, \qquad d_\ell = B^2 m \norm{C}^2  b_{\ell}, \label{eq:bigone_candrhoandbandd}
\end{align}
\hrulefill
\end{figure*}

Theorem~\ref{thm:1} generalizes standard linear convergence results using asynchronous projected gradient descent on functions that satisfy the error bound condition in~\cite{tseng1991rate} to time-invariant feedback optimization problems that satisfy the same error bound condition. 

Theorem~\ref{thm:1} establishes
that linear convergence is attained for static feedback optimization problems,
though the constants and overall convergence rate are different from 
those in the static case. 
That is, linear convergence carries over from conventional optimization
to static feedback optimization. 
To the best of our knowledge this is the first convergence result for time-invariant feedback optimization problems where agents' computations,
communications, and state measurements 
are subject to asynchrony.

Motivated by Theorem~\ref{thm:1}, we next
present a convergence result 
for 
Algorithm~\ref{alg:myAlg} for time-varying feedback optimization problems. 

\begin{theorem} \label{thm:2}
    Let Assumptions~\ref{ass:convex}-\ref{ass:partialAsynch} hold. 
    Fix~$T \in \N$ and fix~${\T = \left\{ \ind{0},\dots, \ind{T} \right\}}$. 
    Suppose that~$\frac{2}{B(L_{x,\ell} + L_{y,\ell}\|C\|^2)} \leq 1$. 
    Then, for all~$\ind{\ell} \in \T$, a stepsize~$\gamma_\ell \in (0,\gamma_{\max,\ell})$, 
    and
    iterates~$\x{}{}{k}$ and~$\y{}{}{k}$ as defined in~\eqref{eq:85} and~\eqref{eq:86}, the sequence~$\{ \x{}{}{k}, \y{}{}{k} \}_{k \in \N_0}$ generated by~$N$ agents executing Algorithm~\ref{alg:myAlg} satisfies
    \begin{multline}
        \at{\eta_{\ell-1} + r_\ell B}{\ell} \leq {a_{\ell} \pt{r_{\ell} - 1}{\ell}}, 
        \bt{\eta_{\ell-1} + r_\ell B}{}     \leq b_{\ell} \pt{r_{\ell} - 1}{\ell}, \\
        \dt{\eta_{\ell-1} + r_\ell B}{}     \leq d_{\ell} \pt{r_{\ell}-1}{\ell}, 
    \end{multline}
where~$\eta_{\ell}$ is from~\eqref{eq:etadef},~$D_{\ell}$ and~$E_{\ell}$ are from~$\eqref{eq:bigone_dande}$,~$F_{\ell}$ is from~\eqref{eq:bigone_f},~$G_{\ell}$ is from~\eqref{eq:bigone_g},~$a_{\ell}$ is from~\eqref{eq:bigone_a},
$\gamma_{\max, \ell}$ is from~\eqref{eq:bigone_gmax}, 
and~$c_{\ell}$,~$\rho_{\ell}$,~$b_{\ell}$, and~$d_{\ell}$ are from~\eqref{eq:bigone_candrhoandbandd}. 
\end{theorem}
\emph{Proof:} See Appendix~C. \hfill $\blacksquare$

Theorem~\ref{thm:2} shows
that the iterates of Algorithm~\ref{alg:myAlg} track solutions
to time-varying feedback optimization problems to within
an error ball. 
We bound errors at time~$\eta_{\ell-1} + r_{\ell}B$ because
this is the last tick of agents' iteration counter~$k$ 
before~$t_{\ell}$ increments to~$t_{\ell+1}$.
Thus, for each~$t_{\ell} \in \T$, 
Theorem~\ref{thm:2} bounds errors
after agents have made all of the progress that
they are going to make towards~$\big(x^*(t_{\ell}), y^*(t_{\ell})\big)$. 

\begin{remark}
Although the results in Theorems~\ref{thm:1} and~\ref{thm:2}
bear some superficial resemblance, the meanings of
these results
differ substantially. 
Specifically, Theorem~\ref{thm:2} bounds~$\alpha$ using~$a_{\ell}$,
which is defined in~\eqref{eq:bigone_a}, where it
can be seen that~$a_{\ell}$ depends on~$a_{\ell-1}$.
Through this recursive dependence, we see that~$a_{\ell}$ depends on~$a_{\ell-1}, a_{\ell-2}, \ldots, a_0$.
Similarly, the value of~$a_{\ell}$ depends
on~$\rho_{\ell-1}, \ldots, \rho_0$, which encode convergence
rates attained for past objective functions, 
and also depends on~$r_{\ell-1}, \ldots, r_0$, which account
for the numbers of agents' computations, communications, and sensor
measurements executed for past objectives. The remaining terms
in the definition of~$a_{\ell}$ in~\eqref{eq:bigone_a} account
for how~$J(\cdot, \cdot; t_{\ell})$ has changed over time.
Thus, Theorem~\ref{thm:2} presents a bound that must quantify
the effects of the entire time history of 
Algorithm~\ref{alg:myAlg} so far, including all 
operations completed by all agents, and the timing
and magnitude of changes in agents' objective function.
Conversely, Theorem~\ref{thm:1} is about a single static
objective and hence has no dependence on the past. 
\hfill $\blacklozenge$
\end{remark}

\begin{remark}
 Theorem~\ref{thm:2} addresses an open problem identified
 in the literature, namely the problem in
 Remark~3 in~\cite{bernstein2018asynchronous}, because,
 in the setting of feedback optimization, 
 we account for agents exhibiting different computational capabilities by allowing asynchronous computations, communications,
 and output measurements.
 That is, agents with different abilities to compute, communicate, and sense may execute these operations at different rates and different 
 times, and Theorem~\ref{thm:2} analyzes their convergence. 
 \hfill $\blacklozenge$ 
\end{remark}

\begin{remark} \label{rem:rell} 
The bounds in Theorem~\ref{thm:2} show that if~$r_{\ell} = 1$ for all~$\ell$,
then agents may not make any progress towards a solution in the sense
that~$\alpha$,~$\beta$, and~$\delta$ may not get smaller. 
This occurs 
because agents' objective is time-varying. Specifically, due to delays
in agents' communications, as agents minimize one cost, say~$J(\cdot, \cdot; t_{\ell})$,
agent~$i$ may send~$x^i_i(k)$ and~$y^i_i(k)$ to other agents. Before that
message is received by other agents, the cost may change to~$J(\cdot, \cdot; t_{\ell+1})$.
Then, if agent~$i$'s message is received by other agents after this change in cost,
those agents receive iterates that do not help it minimize~$J(\cdot, \cdot; t_{\ell+1})$.
It can take up to~$B$ timesteps for all outdated iterates to be received,
and the value~$r_{\ell} = 1$ means that agents only minimize~$J(\cdot, \cdot; t_{\ell})$
for~$B$ timesteps total. It is therefore possible to communicate only
older iterates during that time, and thus
the value~$r_{\ell} = 1$ means that progress towards a minimizer is not guaranteed. 
\hfill $\blacklozenge$ 
\end{remark}

Theorem~\ref{thm:2} analyzes convergence for
agents operating over a fixed time horizon. One may also
ask what performance to expect if Algorithm~\ref{alg:myAlg}
were run indefinitely over an unbounded time horizon. 
To answer this, 
we next present an asymptotic tracking error result for infinite horizon time-varying feedback optimization problems,
and, in light of Remark~\ref{rem:rell}, we consider~$r_{\ell} \geq 2$ for this result. 

\begin{theorem} \label{thm:3}
    Let Assumptions~\ref{ass:convex}-\ref{ass:partialAsynch} hold. Let
    \begin{multline}
        V_\ell \coloneqq 2 \Delta L_t + L_{J,\ell} \sigma_\ell \big( 1+ \Vert C \Vert \big) + \big( M_{x,\ell} + M_{y,\ell} \Vert C \Vert \big) B \diam \big( \X \big) \\
        + \frac{8 E_\ell\left(\frac{G_\ell}{F_\ell} + \frac{E_\ell}{D_\ell}\right) F_\ell B^2 \diam \left(\X\right)^{2}}{2 D_\ell} \big( L_{x,\ell} + L_{y,\ell} \Vert C \Vert^2  \big),
        \end{multline}
        \begin{equation}
        V_\infty \coloneqq \max \bigg\{ a_0, \underset{\ell \geq 0}{\sup} \ V_\ell \bigg\}, \qquad
        \rho_\infty \coloneqq \underset{\ell \geq 0}{\sup} \ \rho_\ell \in (0,1).
    \end{equation}
    Then, the sequence~$\{ \x{}{}{k}, \y{}{}{k} \}_{k \in \N_0}$ generated by~$N$ agents executing Algorithm~\ref{alg:myAlg} with~${r_\ell \geq 2}$ for all~$\ind{\ell} \geq 0$ satisfies
    \begin{align} \label{eq:t3bound}
        \underset{\ell \rightarrow \infty}{\limsup} \ \at{\eta_\ell}{\ell} \leq V_\infty \frac{\pt{}{\infty}}{1-\pt{}{\infty}},
    \end{align}
    where~$\eta_{\ell}$ is from~\eqref{eq:etadef}. 
\end{theorem}
\emph{Proof:} For any~$\eta_\ell$ where~$\ell \geq 0$, Theorem~\ref{thm:2} gives us
the upper bound
    $\at{\eta_\ell}{\ell} \leq a_0 \prod^\ell_{\theta=0} \pt{r_\theta - 1}{\theta} + \sum^\ell_{\tau=1} V_\tau \prod^\ell_{\psi=\tau} \pt{r_\psi - 1}{\psi}$.
Then, using~$V_\tau \leq V_\infty$,~$\pt{}{\tau} \leq \pt{}{\infty}$ for all~$\tau \geq 0$, 
and the fact that~$r_\ell \geq 2$ for all~$\ind{\ell} \in \N$ results in
    ${\at{\eta_\ell}{\ell} \leq V_\infty \bigg( \prod^\ell_{\theta=0} \pt{}{\infty} + \sum^\ell_{\tau=1} \pt{\tau}{\infty} \bigg)}$. 
For~$\ell \rightarrow \infty$, using~$\sum^\infty_{k=1} \nu^k = \frac{\nu}{1-\nu}$ for~$\vert \nu \vert < 1$ gives the result. \hfill $\blacksquare$

Theorem~\ref{thm:3} shows that agents' long-run performance
can be bounded using certain worst-case constants, 
namely~$V_{\infty}$ and~$\rho_{\infty}$.
In particular, 
suppose for some fixed~$t_{\ell}$ that agents
make little to no progress towards
the minimizer of~$J(\cdot, \cdot; t_{\ell})$.
This can be captured mathematically through
larger values of~$V_{\infty}$ and~$\rho_{\infty}$. 
Then Theorem~\ref{thm:3} reveals that 
this lack of progress for a single objective 
may negatively impact long-term performance
by making~$V_{\infty}$ and~$\rho_{\infty}$ larger.
Conversely, consistently high performance (in the sense
of agents making significant progress towards
the minimizer of~$J(\cdot, \cdot; t_{\ell})$ for each~$t_{\ell}$)
produces long-term high performance in a predictable, quantifiable
way by ensuring that both~$V_{\infty}$ and~$\rho_{\infty}$ are smaller. 

\begin{remark}
Bounds similar to~\eqref{eq:t3bound} appear  in~\cite{popkov2005gradient,bernstein2018asynchronous} for certain centralized and decentralized discrete-time correction-only time-varying 
optimization methods. Thus, 
Theorem~\ref{thm:3} shows that we can obtain similar asymptotic 
tracking performance when considering a more general problem formulation that allows asynchrony in agents' computations, communications,
and sensor measurements.  \hfill $\blacklozenge$ 
\end{remark}

In some cases, networks are able to complete their operations with predictable timing. For example, agents may complete their operations at certain frequencies dependent upon onboard computation and communication hardware. In these situations, it can be of interest to design or select certain
specifications, e.g., computational hardware or a network topology, such that the network yields some desired performance. 
For such cases, we next 
provide bounds on the number of operations agents must execute
to attain a certain cost. 
%
\begin{theorem} \label{thm:4}
    Let Assumptions~\ref{ass:convex}-\ref{ass:partialAsynch} hold and let~$\phi > 0$. 
    Fix~$T \in \N$ and fix~$\T = \{t_0, \ldots, t_T\}$. 
    Suppose~$N$ agents are executing Algorithm~\ref{alg:myAlg} with~$r_\ell \equiv r$ for all~$\ind{\ell} \in \T$ and~$r \in \N$ with $r \geq 2$. Let
        $V_{\max} = \max \big\{ a_0, \underset{\ind{\ell} \in \T}{\max} V_\ell \big\}$
        and
        $\rho_{\max} \coloneqq \underset{\ind{\ell} \in \T}{\max} \ \rho_\ell \in (0,1)$.
    If 
    \begin{align}
        r \geq 1 + 
        \frac{ \ln \Big( \frac{V_{\max} \rho_{\max}^{(T+2)(r-1)} + \phi }{V_{\max} +\phi} \Big)}{\ln (\rho_{\max})},
    \end{align}
    then~$\at{\eta_\ell}{\ell} \leq \phi$ for all~$\ind{\ell} \in \T$, where~$\eta_{\ell}$ is from~\eqref{eq:etadef}. 
\end{theorem}
\emph{Proof:} For any~$\eta_\ell$ where~$\ell \geq 0$, Theorem~\ref{thm:2} gives us
the bound
    $\at{\eta_\ell}{\ell} \leq a_0 \prod^\ell_{\theta=0} \pt{r_\theta - 1}{\theta} + \sum^\ell_{\tau=1} V_\tau \prod^\ell_{\psi=\tau} \pt{r_\psi - 1}{\psi}$.
Setting~$r_{\tau}=r$ for all~$\tau \geq 0$, and using the facts~$V_\tau \leq V_{\max}$ and~$\pt{}{\tau} \leq \pt{}{\max}$ for all~$\tau \geq 0$ we can derive the upper bound
    $\at{\eta_\ell}{\ell} \leq V_{\max} \pt{r-1}{\max} \sum^\ell_{\tau=0} \pt{\tau(r-1)}{\max}$.
To upper bound~$\at{\eta_\ell}{\ell}$ by~$\phi$ it is sufficient to have
    $V_{\max} \pt{r-1}{\max} \sum^\ell_{\tau=0} \pt{\tau(r-1)}{\max} \leq \phi$.
From the fact that~$\sum^{n-1}_{m=0} aq^{mc}=a\frac{q^{cn}-1}{q^c - 1}$ when~$q < 1$, 
boundedness by~$\phi$ holds if and only if
    $V_{\max} \pt{r-1}{\max} \frac{\pt{(\ell+1)(r-1)}{\max}-1}{\pt{r-1}{\max} - 1} \leq \phi$.
Setting~$\ell = T$ maximizes the left-hand side. Solving for~$r$ 
completes the proof. 
\hfill $\blacksquare$

Theorem~\ref{thm:2} quantifies agents' performance based on how many operations they complete,
and Theorem~\ref{thm:4} essentially inverts this analysis to determine how many operations
agents must execute to attain a desired level of performance. Theorem~\ref{thm:4}
is for the finite-horizon case, and next we extend it to the asymptotic case.

\begin{corollary} \label{cor:1}
    Let Assumptions~\ref{ass:convex}-\ref{ass:partialAsynch} hold and let~$\phi > 0$. Suppose~$N$ agents execute Algorithm~\ref{alg:myAlg} with~$r_\ell \equiv r$ 
    for all~$\ind{\ell} \geq 0$ and~$r \in \N$ with $r \geq 2$.
    Consider the constants~$V_{\infty}$ and~$\rho_{\infty}$ from Theorem~\ref{thm:3}. 
    If
        $r \geq 1 + \frac{ \ln \Big( \frac{\phi}{V_\infty + \phi} \Big)}{\ln (\rho_{\infty})}$,
    then~$\at{\eta_\ell}{\ell} \leq \phi$ for all~$\ind{\ell} \geq 0$.
\end{corollary}
\emph{Proof:} In Theorem~\ref{thm:3}, use~$\lim_{T \to \infty} q_{\infty}^{\left(T+2\right)c} = 0$. \hfill $\blacksquare$

Both Theorem~\ref{thm:4} Corollary~\ref{cor:1} relate
agents' desired performance bound (in~$\phi$), 
the time-varying nature of the problem (through~$V_{\ell}$), and
agents' rate of convergence (through~$\rho_{\ell}$) to quantify how
the rate at which they complete their operations
(codified in~$r$) affects long-term performance. We next explore
related questions numerically. 

\section{Numerical Results}
\label{sec:simulation}
In this section we consider two examples using Algorithm~\ref{alg:myAlg}. 
The first is a time-varying quadratic program utilizing feedback, and
the second is a network of aircraft modeled with F-16XL dynamics.
\subsection{Time-Varying Quadratic Program with Feedback} \label{subsec:QP}
\begin{figure}[!htb]
    \centering
    \includegraphics[width = 0.40 \textwidth]{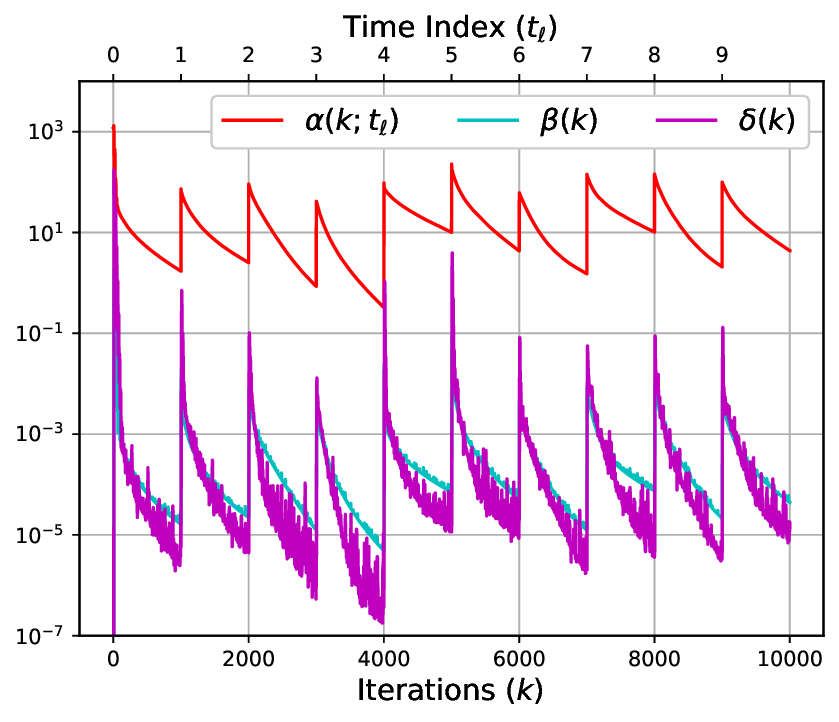}
    \caption{Plot of~$\at{k}{\ell}$,~$\bt{k}{}$, and~$\dt{k}{}$ as a function of~$k$ 
    and~$\ind{\ell} \in \T$ for~${N=10}$ agents executing Algorithm~\ref{alg:myAlg} 
    with~$B=5$. 
    We see that~$\alpha$, $\beta$, and~$\delta$ abruptly increase when agents'
    objective function changes, then gradually decrease between these changes
    as agents complete more computations, communications,
    and sensor readings, which agrees with Theorem~\ref{thm:2}. 
    }
    \label{fig:convergence}
\end{figure}

\begin{figure}[!htb]
    \centering
    \includegraphics[width = 0.44 \textwidth]{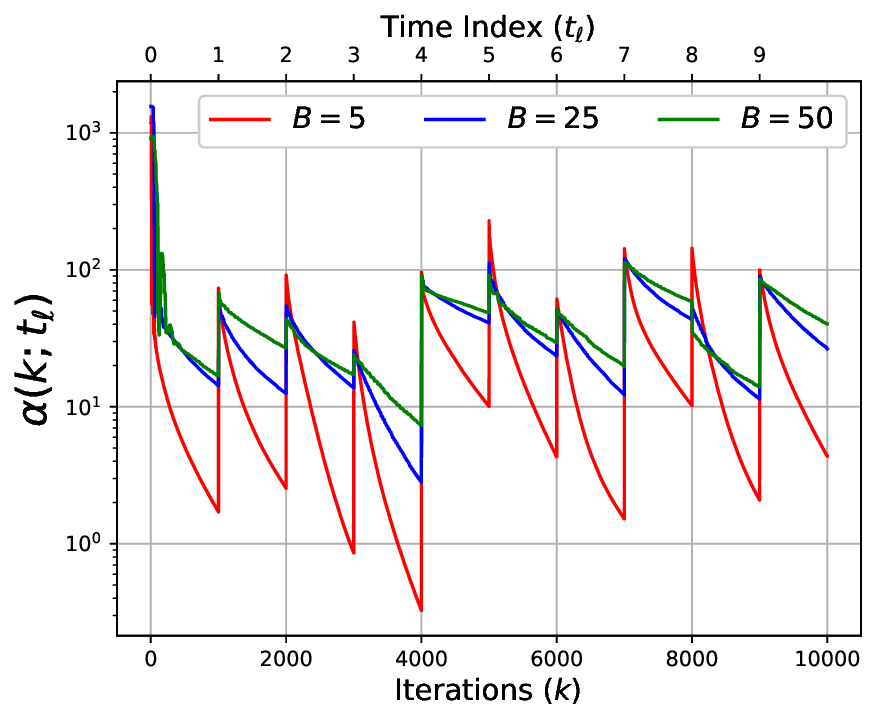}
    \caption{A plot of~$\at{k}{\ell}$ as a function of~$k$ 
    and~$\ind{\ell} \in \T$ for~${N=10}$ agents executing Algorithm~\ref{alg:myAlg} 
    with 
    differing values of~$B$. As~$B$ shrinks, Algorithm~\ref{alg:myAlg} 
    more closely approximates a synchronous algorithm, and convergence
    to a minimizer occurs faster, which
    agrees with Theorem~\ref{thm:2}. 
    }
    \label{fig:multipleB}
\end{figure}

We consider~$N=10$ agents executing Algorithm~\ref{alg:myAlg} as the objective changes. The problem takes the form
\begin{align}
        \underset{x\in \X}{\textnormal{minimize}} & \ \frac{1}{2} x^T Q(\ind\ell) x + q(\ind\ell)^T x + \frac{1}{2} y^T P(\ind\ell) y  + p(\ind\ell)^T y\\
        \textnormal{subject to} & \ y=Cx, 
\end{align}
where~${Q(\ind\ell) \succ 0 \in \R^{20 \times 20}}$,~${P(\ind\ell) \succ 0 \in \R^{10 \times 10}}$, ${q(\ind\ell) \in \R^{20}}$, ${p(\ind\ell) \in \R^{10}}$, and ~${C \in \R^{10 \times 20}}$ are randomly generated. We use the 
set~$\X = [-10,10]^{20}$. 
The objective function changes every~$1,000$ iterations of~$k$.
Each agent is assigned to update two entries of~$x$, 
measure one entry of~$y$, and communicate with other agents.
At each~$k$, agent~$i$ computes an update with probability~$p_{i,u}(k) = 0.01$, measures its output with probability~$p_{i,m}(k) = 0.01$, and communicates with probability~$p_{i,c}(k) = 0.01$.
The stepsize~$\gamma_\ell = 0.001$ was used for all~$\ind{\ell} \in \T$, and we set~$B = 5$. 
In Figures~\ref{fig:convergence} and \ref{fig:multipleB}, Algorithm~\ref{alg:myAlg}
operates on the ``Iterations~$(k)$'' 
timescale (bottom axes), while the objective function 
changes on the ``Time Index~$(\ind{\ell})$'' timescale (top axes).
In Figure~\ref{fig:convergence}, 
the increases in error are due to the changes in the objective function. 
Moreover, we observe that the values of~$\at{k}{\ell},\bt{k}{}$, and~$\dt{k}{}$ trend toward zero between changes in the objective function
for all~$\ind{\ell} \in \T$, which agrees
with the result of Theorem~\ref{thm:2}. 

We demonstrate the effects of the maximum delay length~$B$ in Figure~\ref{fig:multipleB}. As the value of~$B$ decreases, the agents are able to track the minimizers more closely at each time index~$\ind{\ell}$. Intuitively, 
as~$B$ shrinks, we expect Algorithm~\ref{alg:myAlg} to approach a synchronous
feedback optimization algorithm, from which
we expect better tracking performance, 
which we indeed see in Figure~\ref{fig:multipleB}.

\subsection{Aircraft Altitude Tracking}
In this section, Algorithm~\ref{alg:myAlg} is implemented on a network of~$N=8$ F16-XL aircraft to: (1) track a time-varying desired altitude denoted by~$\Phi(\ind{\ell}) \in \R$, 
(2) maintain a specified altitude separation~$\omega_i \in \R$ for all~${i \in [N]}$, and (3) track a time-varying desired acceleration denoted by~$\Psi_i (\ind{\ell}) \in \R$ for all~$i \in [N]$. 
This goal is representative of tracking a target with an \emph{a priori} unknown trajectory while maintaining safe levels of acceleration for the aircraft and avoiding inter-agent collisions. 
The state vector of agent~$i$ is denoted~${x_i = \big[ v_i, \vartheta_i, \varphi_i, \dot{\varphi}_i, \xi_i \big]^T \in \R^5}$ 
for all~$i \in [N]$, where~$\xi_i \in \R$ is the 
altitude,~$v_i \in \R$ is the velocity,~$\vartheta_i \in \R$ is 
the angle of attack,~$\varphi_i \in \R$ is the pitch,
and~$\dot{\varphi}_i \in \R$ is the pitch rate. 
Furthermore, we consider~$y_i = \big[ \dot{v}_i, \xi_i \big]^T$ as 
the outputs for all~$i \in [N]$, where~$\dot{v}_i \in \R$ is the acceleration of the aircraft. 
By linearizing the aircraft longitudinal dynamics~\cite{russell2003non} 
about the operating point~${\bar{x} = \big[ 500 \textnormal{ ft/s}, 0\degree, 0\degree, 0\degree/s, 15,000 \textnormal{ ft} \big]^T \in \R^5}$ we formulated the  input-output relationship~$ y_i = C_i x_i$, where
\begin{align}   C_i &= \begin{bmatrix}
       -0.0133 &-7.3259 &-3.17 &-1.1965 &0.0001 \\
       0& 0& 0& 0& 1
   \end{bmatrix}.
\end{align}
We denote the actual altitude separation 
between aircraft as~${\tilde{\xi}=[\xi_1 - \xi_2, \dots, \xi_i - \xi_{i+1}, \dots, \xi_7 - \xi_8]^T \in \R^7}$ and the desired separation vector~$\omega = [\omega_1,\dots,\omega_7]^T \in \R^7$.
We stack the aircraft states and outputs into the vectors~$x=(x_1^T,\dots,x_8^T)^T \in \R^{40}$ and~$y=(y_1^T,\dots,y_8^T)^T \in \R^{16}$, respectively, to compactly write the network-wide output model as~$y = C x$. Here,~$C\in \R^{16 \times 40}$ is a block diagonal matrix with~$C_i$ on its diagonal for all~$i \in [N]$. Lastly, we stack the desired outputs in a vector~$\Theta (\ind{\ell}) = [\Phi(\ind{\ell}), \Psi_1(\ind{\ell}),\dots, \Phi(\ind{\ell}), \Psi_8(\ind{\ell}) ] \in \R^{16}$.
Agents attempt to satisfy the three goals of this section by tracking the solution of
the following problem: 
\begin{align} \label{prob:planes}
       \!\!\!\!\!\!\! \underset{x\in \X}{\textnormal{minimize}} \ \! \J{x}{y}{\ell} & \!\! \coloneqq \!\!\! \ \frac{1}{2} x^T Q x + \! \frac{1}{2} (y - \Theta(\ind{\ell}))^T \! P (y - \Theta(\ind{\ell}))\\
        & \qquad\qquad\qquad\qquad + \frac{1}{2} (\tilde{\xi} - \omega)^T R (\tilde{\xi} - \omega) \\
        \textnormal{subject to} & \quad \ y = Cx,
\end{align}
where~$Q= 100 \cdot I_{40}\in\R^{40 \times 40}$, $P=  I_8 \otimes \big[ 10^3 \ 0; 0 \ 5\times 10^4 \big]  \in\R^{16 \times 16}$, and $R=10^6 \cdot I_{7} \in\R^{7\times7}$. 
We designed the desired altitude as~${\Phi(\ind{\ell}) = 15000 + 1500 \sin \big( \frac{\ind{\ell}\cdot t_s \cdot \pi}{24} \big) \textnormal{ft}}$,  where~$t_s = 5s$ is the sampling time for the desired altitude, and desired separation as~$\omega_i = 1500 \textnormal{ft}$ for all~$i \in [N]$. 

The desired accelerations for each agent are updated according to~$\Psi_i(\ind{\ell}) = \frac{0.1}{t_s} \big( \Phi(\ind{\ell}) - \frac{1}{N} \sum_{j=1}^N \xi_j^i (\eta_\ell)  \big)$.
In this example, desired outputs are updated every 500 iterations of~$k$ and they are updated $20$ times.
We define~$x_{\max} = [556.2664 \text{ ft/s}, 1.5\degree, 25 \degree, 60 \degree / s,40,000 \text{ ft}]$ and~$x_{\min} = [443.7336 \text{ ft/s}, -13 \degree, -25 \degree, -60 \degree / s, 1000 \text{ ft}]$ so that we have~$\X_i = \{ x \in \R^5 \mid x_{\min,j} \leq x_j \leq x_{\max,j},  j=1,\dots,5 \}$ for all~$i \in [8]$. This yields~$\X = \X_1 \times \dots \times \X_8$.
Agent~$i$ computes an update to~$x^i_i(k)$ with probability~$p_{i,u}(k)=0.5$, measures~$y^i_i(k)$ with probability~$p_{i,m}(k) =0.5$, and communicates its desired state and measured altitude with probability~$p_{i,c}(k) = 0.5$.
The maximum delay is~$B=50$. 

\begin{figure}[!htb]
    \centering
    \includegraphics[width = 0.44 \textwidth]{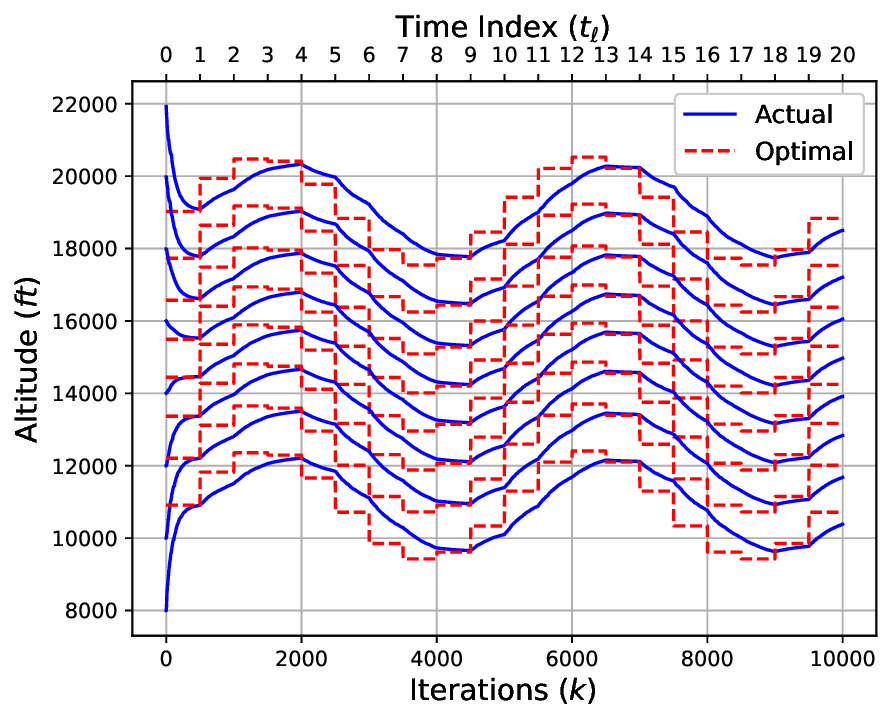}
    \caption{Plot of the optimal agent altitudes (red) and the actual agent altitudes (blue) produced by agents executing Algorithm~\ref{alg:myAlg}. As agents complete more operations their actual altitudes approach the optimal altitudes for all~$\ind{\ell} \in \T$.}
    \label{fig:altitudes}
\end{figure}

\begin{figure}[!htb]
    \centering
    \includegraphics[width = 0.44 \textwidth]{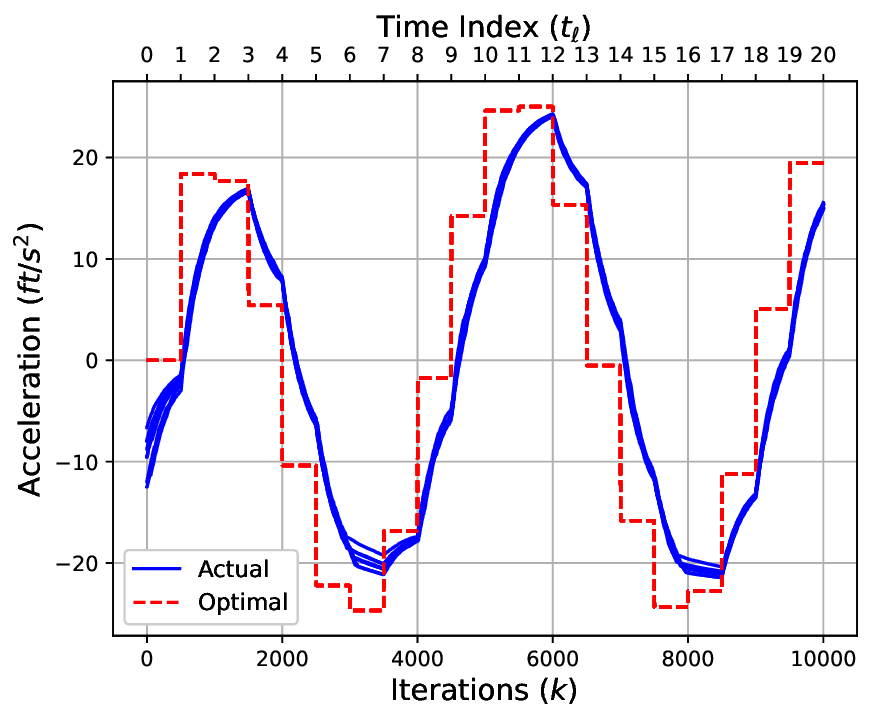}
    \caption{Plot of the optimal agent accelerations (red) and the actual agent accelerations (blue) produced by agents executing Algorithm~\ref{alg:myAlg}. Similar to Figure~\ref{fig:altitudes} agents' accelerations track the optimal accelerations 
    better as agents complete more operations, which is intuitive. 
    }
    \label{fig:accel}
\end{figure}

\begin{figure}[!htb]
    \centering
    \includegraphics[width = 0.44 \textwidth]{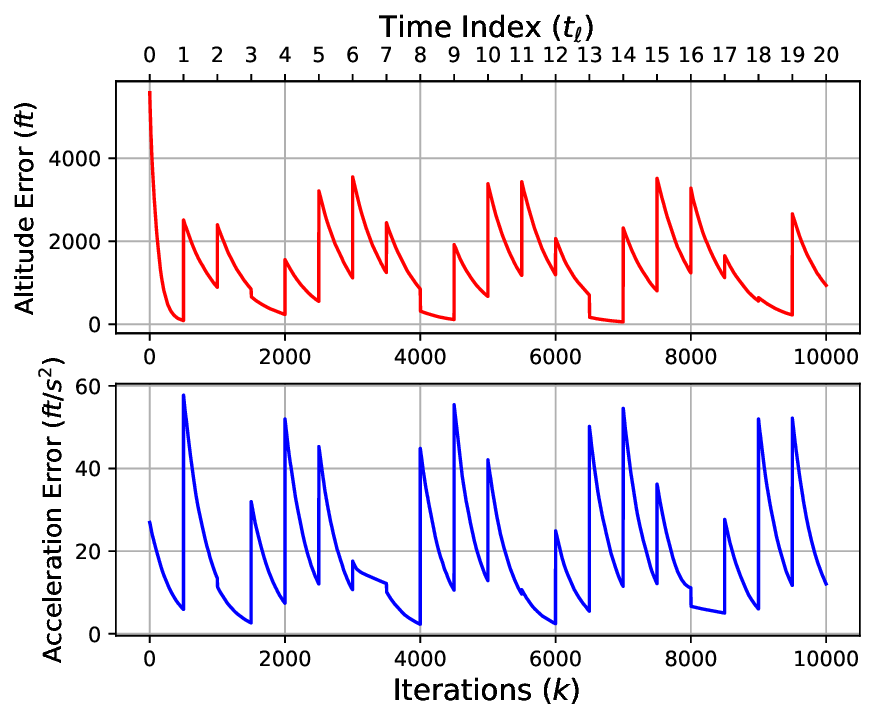}
    \caption{Plot of the error between the optimal outputs and actual outputs. Specifically, at~$k=3999$, the altitude error is equal to~$850.7 ft$ and the acceleration error is~$2.30 ft/s^2$. These values
    are the errors that result both from the time-varying nature of the agents'
    problem and the asynchrony in their operations. 
    }
    \label{fig:error}
\end{figure}

Figure~\ref{fig:altitudes} shows the optimal altitudes of each agent along with the actual agent altitudes as agents track the solution of~\eqref{prob:planes}. Note that the optimal altitudes are not equal to the desired altitude~$\Psi(\ind{\ell})$ because we include the altitude separation term~$\omega$ in the objective function. Even under a high degree of asynchrony (i.e., large~$B$)
we observe that the agents are able to track the optimal altitudes closely.
Figure~\ref{fig:accel} shows the agents' actual accelerations alongside the optimal accelerations. As with altitudes, agents are able to track the optimal accelerations,
even under asynchrony. 
We show the errors between the optimal and actual altitudes and accelerations in Figure~\ref{fig:error}. 

\section{Conclusion}
\label{sec:conclusion}
We presented a decentralized asynchronous algorithm for tracking the solution of 
time-varying
feedback optimization problems, and we derived tracking error bounds for a
class of problems. 
This work extended the feedback optimization concept to distributed problems
with asynchrony, and it showed the viability of feedback optimization
as a means to provide robustness to disturbances that result from asynchrony. 
Future work will explore
analogous results for problems with more complex
observation maps and problems in which agents asynchronously
sample an objective function. 
\appendix
\begin{appendices}
\subsection{Technical Lemmas}
This section provides lemmas used
in Appendices~\ref{sec:t1proof} and~\ref{sec:t2proof}. 
\begin{lemma}[Descent Lemma]{\cite[Proposition A.32]{bertsekas1989parallel}} 
\label{lem:descent}
If~${f:\R^n \rightarrow \R}$ is continuously differentiable and has the property~${\Vert \nabla f(x) - \nabla f(y) \Vert \leq K \Vert x-y \Vert}$ for every~${x,y \in \R^n}$, then
                $f\left(x+y\right) \leq f\left(x\right)+y^{T}\nabla f\left(x\right)+\frac{K}{2}\left\Vert y\right\Vert^{2}$. \hfill $\blacksquare$
\end{lemma}
\begin{lemma} \label{lem:projection}
For all~$\ind{\ell} \in \T$, all~$i \in [N]$, and all~$k \geq 0$ we have the bound
     $\s{i}{k}{\ell}^T \nabla_{\x{}{i}{}} \J{\x{i}{}{k}}{\y{i}{}{k}}{\ell} \leq  -\frac{1}{\gamma_\ell} \norm{\s{i}{k}{\ell}}^2$.
\end{lemma}
\emph{Proof:} For~$v \in \X \subset \R^n$, $x \in \R^n \backslash \X$, and~$z=\Pi_\X \left[x\right]$, a property of orthogonal 
projections~\cite{cheney1959proximity} is~$0 \geq (z-x)^T(z-v)$.
Define~${z \coloneqq \Pi_{\X_i}\big[ \x{i}{i}{k} - \gamma_\ell \nabla_{\x{}{i}{}} \J{\x{i}{}{k}}{\y{i}{}{k}}{\ell}\big]}$. Then
we find that
    $0 \geq \big(z - \big( \x{i}{i}{k} - \gamma_\ell \nabla_{\x{}{i}{}} \J{\x{i}{}{k}}{\y{i}{}{k}}{\ell} \big) \big)^T\big(z-\x{i}{i}{k} \big)$.
Expanding the inner product and combining like terms then gives us
    $0 \geq \norm{z-\x{i}{i}{k}}^2 
     + \big( z - \x{i}{i}{k} \big)^T \gamma_\ell \nabla_{\x{}{i}{}} \J{\x{i}{}{k}}{\y{i}{}{k}}{\ell}$.  
Using~$\s{i}{k}{\ell} = z - \x{i}{i}{k}$
and rearranging completes the proof. \hfill $\blacksquare$
\begin{lemma}[\!\!{\cite[Section 7.5.1]{bertsekas1989parallel}}] \label{lem:outOfDate}
For all~$\ind{\ell} \in \T$, all~$i \in [N]$, and all~$k \geq 0$, we have
           $\norm{ \x{i}{}{k} - \x{}{}{k} } \leq \sum^{k-1}_{\tau = k - B} \norm{\s{}{\tau}{\ell}}$. \hfill $\blacksquare$
\end{lemma}

\begin{lemma} \label{lem:measurement}
For all~$\ind{\ell} \in \T$, all~$i \in [N]$, and all~$k \geq 0$, we have
    $\norm{\y{i}{}{k} - \y{}{}{k}} \leq N \norm{C} \sum^{k-1}_{\tau = k-B} \norm{\s{}{\tau}{\ell}}$.
\end{lemma}

\emph{Proof:} We define~$\s{}{\tau}{0}=0$ for~$\tau<0$. For any~${i\in [N]}$, 
the triangle inequality gives
    $\norm{\y{}{}{k} - \y{i}{}{k}}  \leq \sum^N_{j=1} \Vert \y{}{j}{k} - \y{j}{j}{\mu^i_j(k)} \Vert$.
Using ${\y{}{j}{k} = C_{j*} \x{}{}{k}}$ and~$\y{j}{j}{\mu^i_j(k)}=C_{j*} \x{}{}{\mu^i_j(k)}$, where  ${C_{j*} \in \R^{m_j \times n}}$, we then
find the upper bound 
    ${\norm{\y{}{}{k} - \y{i}{}{k}} \leq \sum^N_{j=1} \big\Vert C_{j*} \x{}{}{k} - C_{j*} \x{}{}{\mu^i_j(k)} \big\Vert}$.
From $\norm{AB} \leq \norm{A} \norm{B}$ and the fact that~$\norm{C_{j*}} \leq \norm{C}$ for all~$j$, we have
    $\norm{\y{}{}{k} - \y{i}{}{k}} \leq \norm{C} \sum^N_{j=1} \Vert \sum^{k-1}_{\tau = \mu^i_j(k)} \s{}{\tau}{\ell} \Vert$.
By partial asynchrony~$k-B \leq \mu^i_j(k)$, and using this and the triangle inequality gives
    $\norm{\y{}{}{k} - \y{i}{}{k}}  \leq \norm{C} \sum^N_{j=1} \sum^{k-1}_{\tau = k-B} \norm{\s{}{\tau}{\ell}}$. \hfill $\blacksquare$
\begin{lemma} \label{lem:cost}
For all~$\ind{\ell} \in \T$ and~$k \in \{\eta_{\ell-1}, \ldots, \eta_{\ell}-B\}$, we have
\begin{multline} 
    \J{\x{}{}{k+B}}{\y{}{}{k+B}}{\ell} - \J{\x{}{}{k}}{\y{}{}{k}}{\ell} \leq  \\  
     - \left( \frac{  2 - \gamma_{\ell} \big((1 + B) \Lx + (1 + B N) \norm{C}^2 \Ly \big) }{2 \gamma_{\ell}} \right) 
     \!\sum^{k+B-1}_{\tau = k} \! \big\Vert \s{}{\tau}{\ell} \big\Vert^2 
     \\ +  N B \frac{\Lx+\Ly N \norm{C}^2}{2} \sum^{k-1}_{\tau = k - B} \big\Vert \s{}{\tau}{\ell} \big\Vert^2.
\end{multline}
\end{lemma}

\emph{Proof:} We first bound~$\f{\x{}{}{k+1}}{\ell} - \f{\x{}{}{k}}{\ell}$. 
By definition,~$\f{\x{}{}{k+1}}{\ell} = \f{\x{}{}{k} + \s{}{k}{\ell}}{\ell}$.
Then, by Lemma~\ref{lem:descent},
    ${\f{\x{}{}{k\!+\!1}}{\ell} \!\leq\! \f{\x{}{}{k}}{\ell} 
    \!+\! \s{}{k}{\ell}^T \nabla_{\x{}{}{}} \f{\x{}{}{k}}{\ell} 
    \!+\! \frac{\Lx}{2} \norm{\s{}{k}{\ell}}^2}$\!\!.
Expressing the inner product as a sum, 
adding~$\nabla_{\x{}{i}{}} \f{\x{i}{}{k}}{\ell} - \nabla_{\x{}{i}{}} \f{\x{i}{}{k}}{\ell}$ inside
the sum, and rearranging, we find
\begin{multline}
    \f{\x{}{}{k+1}}{\ell}  \leq \f{\x{}{}{k}}{\ell}
    + \sum^N_{i=1}\s{i}{k}{\ell}^T  \nabla_{\x{}{i}{}} \f{\x{i}{}{k}}{\ell} \\    
    + \sum^N_{i=1}\s{i}{k}{\ell}^T \big(\nabla_{\x{}{i}{}} \f{\x{}{}{k}}{\ell} - \nabla_{\x{}{i}{}} \f{\x{i}{}{k}}{\ell}\big) 
    + \frac{\Lx}{2} \norm{\s{}{k}{\ell}}^2.
\end{multline}
In the second sum, we apply the Cauchy-Schwarz inequality, 
the Lipschitz property of~$\nabla_{\x{}{}{}} \f{\cdot}{\ell}$, and 
the result of Lemma~\ref{lem:outOfDate} to find
\begin{multline}
    \f{\x{}{}{k+1}}{\ell}  \leq \f{\x{}{}{k}}{\ell}  
    + \sum^N_{i=1}\s{i}{k}{\ell}^T  \nabla_{\x{}{i}{}} \f{\x{i}{}{k}}{\ell} \\    
    + \Lx \sum^N_{i=1} \big\Vert \s{i}{k}{\ell} \big\Vert \sum^{k-1}_{\tau = k - B} \big\Vert \s{}{\tau}{\ell} \big\Vert 
    + \frac{\Lx}{2} \norm{\s{}{k}{\ell}}^2.
\end{multline}
Using~${\big\Vert \s{i}{k}{\ell} \big\Vert \cdot \big\Vert \s{}{\tau}{\ell} \big\Vert  \leq  \frac{1}{2} \big( \big\Vert \s{i}{k}{\ell} \big\Vert^2 + \big\Vert \s{}{\tau}{\ell} \big\Vert^2} \big)$
gives
%
%
\begin{multline} \label{eq:fDescent}
    \f{\x{}{}{k+1}}{\ell}  \leq \f{\x{}{}{k}}{\ell} + \sum^N_{i=1}\s{i}{k}{\ell}^T  \nabla_{\x{}{i}{}} \f{\x{i}{}{k}}{\ell} \\
    \!\!\!\!+ \frac{\Lx}{2} \Bigg( B \norm{\s{}{k}{\ell}}^2 
    + N \sum^{k-1}_{\tau = k - B} \big\Vert \s{}{\tau}{\ell} \big\Vert^2 \Bigg) 
    + \frac{\Lx}{2} \norm{\s{}{k}{\ell}}^2.
\end{multline}

Next, we can bound~$\g{\y{}{}{k+1}}{\ell} - \g{\y{}{}{k}}{\ell}$ 
using the same steps with Lemma~\ref{lem:measurement} in place of Lemma~\ref{lem:outOfDate}. 
This leads to
\begin{multline} \label{eq:gDescent}
    \g{\y{}{}{k+1}}{\ell} \leq  \g{\y{}{}{k}}{\ell} 
    + \sum^N_{i=1} \s{i}{k}{\ell} ^T \Ci^T \nabla_{\y{}{}{}} \g{\y{i}{}{k}}{\ell} \\
    + \frac{\Ly N \norm{C}^2}{2} \Bigg( \! B \|s(k)\|^2 
    + N \!\!\!\! \sum^{k-1}_{\tau = k - B} \!\!\! \Vert s(\tau) \Vert^2 \! \Bigg)     
    + \frac{\Ly \norm{C}^2}{2} \norm{s(k)}^2.
\end{multline}
Adding~\eqref{eq:fDescent} and~\eqref{eq:gDescent}
and applying Lemma~\ref{lem:projection} gives
\begin{multline} 
    \J{\x{}{}{k+1}}{\y{}{}{k+1}}{\ell} - \J{\x{}{}{k}}{\y{}{}{k}}{\ell} \leq \\    
    \Bigg( {-\frac{1}{\gamma_{\ell}}} + \frac{\Lx + \Ly \norm{C}^2}{2}  + B\frac{ \Lx + \Ly N \norm{C}^2}{2} \Bigg)  \norm{\s{}{k}{\ell}}^2 \\    
    +  N \frac{\Lx+\Ly N \norm{C}^2}{2} \sum^{k-1}_{\tau = k - B} \big\Vert \s{}{\tau}{\ell} \big\Vert^2. 
\end{multline}
We apply this to~$k, k+1, \dots, k+B - 1$ and sum to conclude. 
\hfill $\blacksquare$
\begin{lemma} \label{lem:update}
Define~$\chi(\gamma_{\ell}) = 2\gamma_\ell \big( 1+\gamma_\ell \big) \Lt^2 B N$.
For each~$\ind{\ell} \in \T$ and all~$k \in \{\eta_{\ell-1}, \ldots, \eta_{\ell}-B\}$ we have 
\begin{multline}
    \big\Vert \x{}{}{k+1} - \x{}{}{k} \big\Vert^2  
    \leq 
    \frac{ 1 + \chi(\gamma_{\ell})}{1-\gamma_\ell}  \sum_{\tau=k}^{k+B-1} \Vert \s{}{\tau}{\ell} \Vert^2 \\    
    + \frac{\chi(\gamma_{\ell})}{1-\gamma_\ell}  
    \left(\sum_{\tau=k}^{k+B-1} \Vert \q{}{\tau} \Vert^2 + \!\!\sum_{\tau=k-B}^{k-1} \Big(\Vert \s{}{\tau}{\ell} \Vert^2
    + \Vert \q{}{\tau} \Vert^2\Big) \right). 
\end{multline}
\end{lemma}

\emph{Proof:} Fix some $k\in\left\{ \eta_{\ell-1},\dots,\eta_{\ell}\right\}$.
For each $i\in[N]$ let $k^{i}$ be the smallest
element of $\K^{i}$ that exceeds $k$. Then, for each~$i \in [N]$,
\begin{align}
    \!\!\!\!\!\!\!\!\!\!\x{i}{i}{k^i} \!&=\! \x{i}{i}{k} \label{eq:1} \\
    \!\!\!\!\!\!\!\!\!\!\s{i}{k^i}{\ell} \!&=\! \Pi_{\X_i} \!\Big[ \x{i}{i}{k^i} \!-\! \gamma_\ell \nabla_{\x{}{i}{}} \J{\x{i}{}{k^i}}{\y{i}{}{k^i}}{\ell}  \Big] \!-\! \x{i}{i}{k^i}. \label{eq:2}
\end{align}
Using the triangle inequality and
the non-expansive property of the projection operator we have
\begin{align}
   &\big\Vert \Pi_{\X_i} \big[ \x{i}{i}{k} - \gamma_\ell \nabla_{\x{}{i}{}} \J{\x{}{}{k}}{\y{}{}{k}}{\ell} \big] - \x{i}{i}{k} \big\Vert 
   \leq \\
    & \gamma_\ell \big\Vert \nabla_{\x{}{i}{}} \J{\x{}{}{k}}{\y{}{}{k}}{\ell} - \nabla_{\x{}{i}{}} \J{\x{i}{}{k^i}}{\y{i}{}{k^i}}{\ell}  \big\Vert \\
   &+  \big\Vert \Pi_{\X_i} \big[ \x{i}{i}{k} - \gamma_\ell \nabla_{\x{}{i}{}} \J{\x{i}{}{k^i}}{\y{i}{}{k^i}}{\ell} \big] - \x{i}{i}{k} \big\Vert.
\end{align}
Subtracting the~$\gamma_\ell$ term from both sides and using~\eqref{eq:1} and~\eqref{eq:2} gives
\begin{align}
    &\big\Vert \s{i}{k^i}{\ell} \big\Vert \geq \big\Vert \Pi_{\X_i} \big[ \x{i}{i}{k} - \gamma_\ell \nabla_{\x{}{i}{}} \J{\x{}{}{k}}{\y{}{}{k}}{\ell} - \x{i}{i}{k} \big] \big\Vert \\
    & - \gamma_\ell \big\Vert \nabla_{\x{}{i}{}} \J{\x{}{}{k}}{\y{}{}{k}}{\ell} - \nabla_{\x{}{i}{}} \J{\x{i}{}{k^i}}{\y{i}{}{k^i}}{\ell}  \big\Vert.
\end{align}
Using the Lipschitz property of~$\nabla_x \J{\cdot}{\cdot}{\ell}$ 
%
and using~$(a-\gamma b)^2 \geq (1-\gamma)a^2 - \gamma (1+\gamma)b^2$ on the right-hand side then gives us
\begin{multline} \label{eq:6}
    \big\Vert \s{i}{k^i}{\ell} \big\Vert^2 \geq - \gamma_\ell \big( 1+\gamma_\ell \big) \Lt^2 \big\Vert \big( \x{}{}{k}, \y{}{}{k}\big) - \big( \x{i}{}{k^i},\y{i}{}{k^i} \big) \big\Vert^2
    \\ \big(1-\gamma_\ell \big)\big\Vert \Pi_{\X_i} \big[ \x{i}{i}{k} 
    - \gamma_\ell \nabla_{\x{}{i}{}} \J{\x{}{}{k}}{\y{}{}{k}}{\ell}  \big] - \x{i}{i}{k} \big\Vert^2.     
\end{multline}
By definition,~$k\leq k^i \leq k+B-1$. By partial asynchronism, we also have~$k-B \leq \tau^i_j(k^i) \leq k+B-1$ and~$k-B \leq \mu^i_j(k^i) \leq k+B-1$.
Then, for all~$j$,
\begin{align}
    \Bigg\Vert \begin{bmatrix}
    \x{j}{j}{k} \\ \y{}{j}{k}
    \end{bmatrix} - \begin{bmatrix}
    \x{j}{j}{\tau^i_j(k^i)} \\ \y{j}{j}{\mu^i_j(k^i)}
    \end{bmatrix} \Bigg\Vert^2 \leq \Bigg( \sum_{\tau=k-B}^{k+B-1} \Bigg\Vert \begin{bmatrix}
        \s{j}{\tau}{\ell} \\
        \q{j}{\tau}
    \end{bmatrix} \Bigg\Vert   \Bigg)^2.
\end{align}
Using the inequality $(a_1+\dots+a_{2B})^2 \leq 2B (a_1^2 + \dots + a_{2B}^2)$,
summing over all~$j \in [N]$, and using~\eqref{eq:3} and~\eqref{eq:4} 
then gives
\begin{align} \label{eq:5}
    \Bigg\Vert \begin{bmatrix}
    \x{}{}{k} \\ \y{}{}{k}
    \end{bmatrix} - \begin{bmatrix}
    \x{i}{}{k^i} \\ \y{i}{}{k^i}
    \end{bmatrix} \Bigg\Vert^2 \leq 2B \sum_{\tau=k-B}^{k+B-1} \Bigg\Vert \begin{bmatrix}
        \s{}{\tau}{\ell} \\
        \q{}{\tau}
    \end{bmatrix} \Bigg\Vert^2.
\end{align}
Applying~\eqref{eq:5} to~\eqref{eq:6} gives
\begin{multline} \label{eq:7}
    \big\Vert \s{i}{k^i}{\ell} \big\Vert^2 \geq - 2\gamma_\ell \big( 1+\gamma_\ell \big) \Lt^2 B \sum_{\tau=k-B}^{k+B-1} \Bigg\Vert \begin{bmatrix}
        \s{}{\tau}{\ell} \\
        \q{}{\tau}
    \end{bmatrix} \Bigg\Vert^2 + \\
    \!\!\!\!\big(1\!-\!\gamma_\ell \big)\big\Vert \Pi_{\X_i} \big[ \x{i}{i}{k} 
    - \gamma_\ell \nabla_{\x{}{i}{}} \!\J{\x{}{}{k}}{\y{}{}{k}}{\ell}  \big] \!-\! \x{i}{i}{k} \big\Vert^2.     
\end{multline}

Next, we sum the left-hand side over~$i \in [N]$ to find
\begin{equation} \label{eq:RHSbound}
    \sum^N_{i=1} \big\Vert \s{i}{k^i}{\ell} \big\Vert^2 \leq \sum^N_{i=1} \sum^{k+B-1}_{\tau=k} \big\Vert \s{i}{\tau}{\ell} \big\Vert^2 = \sum^{k+B-1}_{\tau=k} \big\Vert \s{}{\tau}{\ell} \big\Vert^2.
\end{equation}
Summing over~$i \in [N]$ in~\eqref{eq:7} and applying~\eqref{eq:RHSbound} gives 
\begin{multline}
    \sum^{k+B-1}_{\tau=k} \big\Vert \s{}{\tau}{\ell} \big\Vert^2 \geq \big(1-\gamma_\ell \big)\big\Vert \x{}{}{k+1} - \x{}{}{k} \big\Vert^2 \\
    {-2}\gamma_\ell \big( 1+\gamma_\ell \big) \Lt^2 B N \sum_{\tau=k-B}^{k+B-1} \Bigg\Vert \begin{bmatrix}
        \s{}{\tau}{\ell} \\
        \q{}{\tau}
    \end{bmatrix} \Bigg\Vert^2.
\end{multline}
The result follows by 
taking~${\gamma_\ell \in (0,1)}$, and dividing by~$1-\gamma_\ell$. \hfill $\blacksquare$
\begin{lemma} \label{lem:technical} 
Define~$K_{\ell} = \Lx + \Ly\|C\|^2$. 
For all~$\ind{\ell} \in \T $,~${x \in \R^n}$,~$y \in \R^m$, $\bar{x}\in\X$ with~$\bar{y} = C \bar{x}\in\Y$,
${\x{1}{}{},\dots,\x{N}{}{} \in \R^n}$, and~$\y{1}{}{},\dots,\y{N}{}{} \in \R^m$, we have 
\begin{multline}
     \J{a}{b}{\ell}-\J{\bar{x}}{\bar{y}}{\ell} \leq \\ 
    \bigg( N  \bigg( \frac{3}{2} \Big( K_{\ell} + \frac{1}{\gamma_\ell} \Big)^2 + \frac{5}{2} \bigg)\gamma^2_\ell \Lt^2 + \frac{3}{2} \Lx^2 \bigg) \sum^N_{j=1} \norm{x - \x{j}{}{}}^2 \\ + \bigg( N \bigg( \frac{3}{2} \Big( K_{\ell} + \frac{1}{\gamma_\ell} \Big)^2 + \frac{5}{2} \bigg)\gamma^2_\ell \Lt^2 + \frac{3}{2} \Ly^2 \norm{C}^2 \bigg) \sum^N_{j=1} \norm{ y - \y{j}{}{}}^2    \\
    + N \bigg( \frac{3}{2} \Big( K_{\ell} + \frac{1}{\gamma_\ell} \Big)^2 + \frac{5}{2} \bigg) \norm{\bar{a} - x}^2  
    + N \bigg( \frac{3}{2} \Big( K_{\ell} \Big)^2  +\frac{5}{2} \bigg) \norm{\bar{x} - x}^2,
\end{multline}
where~$\bar{a}\coloneqq \Pi_{\X}\big[x - \gamma_\ell \nabla_x \J{x}{y}{\ell} \big]$ and~$a$ is a vector with components~$a_i \coloneqq \Pi_{\X_i} \big[ \x{}{i}{} - \gamma_\ell \nabla_{\x{}{i}{}} \J{\x{i}{}{}}{\y{i}{}{}}{\ell} \big]$
for all~$i \in [N]$, along with~$\bar{b} = C\bar{a}$ and~$b=Ca$.
\end{lemma}

\emph{Proof:} For each~$i \in [N]$, since we define~$a_i$ to be the orthogonal projection of~$\x{}{i}{} - \gamma_\ell \nabla_{\x{}{i}{}} \J{\x{i}{}{}}{\y{i}{}{}}{\ell}$ onto~$\X_i$ and~$\bar{x}_i \in \X_i$, we have 
    $0 \leq \big\langle a_i - \bar{x}_i, \x{}{i}{} - \gamma_\ell \nabla_{\x{}{i}{}} \J{\x{i}{}{}}{\y{i}{}{}}{\ell} - a_i \big\rangle$,
or, equivalently,
\begin{equation} \label{eq:8}
    \!\!\!\!0 \!\leq\! \big\langle a_i \!-\! \bar{x}_i,  -\nabla_{\x{}{i}{}} \f{\x{i}{}{}}{\ell} \!-\! \Ci^T \nabla_y \g{\y{i}{}{}}{\ell} 
    \!+\! \frac{1}{\gamma_\ell} ( \x{}{i}{} - a_i ) \big\rangle.
\end{equation}

The definition of~$\J{\cdot}{\cdot}{\ell}$ gives 
    $\J{a}{b}{\ell}-\J{\bar{x}}{\bar{y}}{\ell} = \f{a}{\ell} - \f{\bar{x}}{\ell} + \g{b}{\ell} - \g{\bar{y}}{\ell}$.
By the Mean Value Theorem there exists a point~$(\phi,\nu)$ between~$(a,b)$ and~$(\bar{x},\bar{y})$ such that
\begin{equation}
    \J{a}{b}{\ell}-\J{\bar{x}}{\bar{y}}{\ell} = \big\langle a-\bar{x}, \nabla_x \f{\phi}{\ell} \big\rangle 
    + \big\langle b-\bar{y}, \nabla_y \g{\nu}{\ell} \big\rangle.
\end{equation}
Using~$\bar{y}=C\bar{x}$ and~$b=Ca$ and simplifying we have
\begin{equation}
    \!\!\J{a}{b}{\ell} - \J{\bar{x}}{\bar{y}}{\ell} \!=\!\! \sum^N_{i=1} \! \big\langle a_i-\bar{x}_i, \!\nabla_{x_i} \f{\phi}{\ell} + \Ci^T\!\nabla_y \g{\nu}{\ell}\!\big\rangle .
\end{equation}
Using~\eqref{eq:8} we add~$-\Big( \nabla_{\x{}{i}{}} \f{\x{i}{}{}}{\ell} + \Ci^T \nabla_y \g{\y{i}{}{}}{\ell} \Big) + \frac{1}{\gamma_\ell} ( \x{}{i}{} - a_i )$ inside the second argument of the inner product. Then 
\begin{multline}
    \J{a}{b}{\ell}-\J{\bar{x}}{\bar{y}}{\ell} \leq \sum^N_{i=1} \Vert a_i-\bar{x}_i \Vert \cdot \Big\Vert \nabla_{\x{}{i}{}} \f{\phi}{\ell} \\
    + \Ci^T\nabla_y \g{\nu}{\ell} 
    -\nabla_{\x{}{i}{}} \f{\x{i}{}{}}{\ell} - \Ci^T \nabla_y \g{\y{i}{}{}}{\ell}  
    + \frac{1}{\gamma_\ell} ( \x{}{i}{} - a_i ) \Big\Vert.
\end{multline}
Using the triangle inequality,~$\norm{\Ci} \leq \norm{C}$, 
%
the Lipschitz property of~$\nabla_x \f{\cdot}{\ell}$ and~$\nabla_y \g{\cdot}{\ell}$, and 
$\big\Vert \nabla_{\x{}{i}{}} \f{\phi}{\ell} -\nabla_{\x{}{i}{}} \f{\x{i}{}{}}{\ell} \big\Vert \leq \big\Vert \nabla_{x} \f{\phi}{\ell} -\nabla_{x} \f{\x{i}{}{}}{\ell} \big\Vert$ for all~$i \in [N]$, we find 
\begin{multline}
    \J{a}{b}{\ell}-\J{\bar{x}}{\bar{y}}{\ell} \leq \sum^N_{i=1} \Vert a_i-\bar{x}_i \Vert \cdot \Big( \Lx \Vert \phi - \x{i}{}{} \Vert \\
    + \Ly \norm{C} \Vert \nu - \y{i}{}{} \Vert + \frac{1}{\gamma_\ell} \big\Vert \x{}{i}{} - a_i \big\Vert  \Big).
\end{multline}
We have both~${\Vert a_i-\bar{x}_i \Vert \leq \Vert a-\bar{x} \Vert}$ and~${\big\Vert \x{}{i}{} - a_i \big\Vert \leq \big\Vert x - a \big\Vert}$
for all~$i \in [N]$, 
which give us
\begin{multline}
    \J{a}{b}{\ell}-\J{\bar{x}}{\bar{y}}{\ell} \leq \sum^N_{i=1} \Vert a -\bar{x} \Vert \cdot \Big( \Lx \Vert \phi - \x{i}{}{} \Vert \\
    + \Ly \norm{C} \Vert \nu - \y{i}{}{} \Vert + \frac{1}{\gamma_\ell} \big\Vert x - a \big\Vert  \Big). \label{eq:11}
\end{multline}

We continue by bounding terms on the right-hand side of~\eqref{eq:11} separately. Since~$\phi$ 
lies between
the points~$a$ and~$\bar{x}$, we have 
${\Vert \phi - x \Vert \leq \Vert a - x \Vert + \norm{\bar{x} - x}}$ for all~$x \in \R^n$. 
Then 
\begin{align}
    \Vert \phi - \x{i}{}{} \Vert & = \Vert \phi - x + x - \x{i}{}{} \Vert \leq \Vert \phi - x \Vert + \norm{x - \x{i}{}{}} \\
    &\leq \Vert a - x \Vert + \norm{\bar{x} - x} + \norm{x - \x{i}{}{}}. \label{eq:9}
\end{align}
Similar reasoning gives
\begin{align}
    \Vert \nu - \x{i}{}{} \Vert 
    &\leq \Vert b - y \Vert + \Vert \bar{y} - y \Vert + \norm{y - \y{i}{}{}}. \label{eq:10}
\end{align}
Applying~\eqref{eq:9} and~\eqref{eq:10} to~\eqref{eq:11}, 
substituting in~$b=Ca$,~$y=Cx$, and~$\bar{y}=C\bar{x}$, and
combining like terms, we find
\begin{multline}
    \J{a}{b}{\ell}-\J{\bar{x}}{\bar{y}}{\ell} \leq \sum^N_{i=1} \Vert a -\bar{x} \Vert \cdot \\
    \Big( \big( \Lx + \Ly \norm{C}^2 + \frac{1}{\gamma_\ell} \big)  \Vert a - x \Vert 
    + \big( \Lx + \Ly \norm{C}^2 \big) \norm{\bar{x} - x} \\
    + \Lx\norm{x - \x{i}{}{}} 
    + \Ly \norm{C}\norm{y - \y{i}{}{}} 
     \Big).
\end{multline}
By adding zero and using the triangle inequality we reach
\begin{multline}
    \J{a}{b}{\ell}-\J{\bar{x}}{\bar{y}}{\ell} \leq  \sum^N_{i=1} \Big( \norm{a - \bar{a}} + \norm{\bar{a} - x} + \norm{x -\bar{x}} \Big) \\
    \cdot \bigg(( \Lx + \Ly \norm{C}^2) \norm{\bar{x} - x} + \Lx\norm{x - \x{i}{}{}} + \Ly \norm{C}\norm{y - \y{i}{}{}}  + \\
    \Big( \Lx + \Ly \norm{C}^2 + \frac{1}{\gamma_\ell} \Big) ( \norm{a - \bar{a}} + \norm{\bar{a} - x})\bigg). 
\end{multline}
Expanding, using~$ab\leq\frac{1}{2}\left(a^{2}+b^{2}\right)$, and simplifying, 
we get 
\begin{multline}
    \J{a}{b}{\ell}-\J{\bar{x}}{\bar{y}}{\ell} \leq \\ 
    N \bigg( \frac{3}{2} \Big( \Lx + \Ly \norm{C}^2 + \frac{1}{\gamma_\ell} \Big)^2 + \frac{5}{2} \bigg)
    \Big( \norm{a - \bar{a}}^2 + \norm{\bar{a} - x}^2  \Big) \\
    + N \bigg( \frac{3}{2} \Big( \Lx + \Ly \norm{C}^2 \Big)^2 
    + \frac{5}{2} \bigg) \norm{\bar{x} - x}^2 \\
    +  \frac{3}{2} \Lx^2 \sum^N_{i=1} \norm{x - \x{i}{}{}}^2 
    + \frac{3}{2} \Ly^2 \norm{C}^2  \sum^N_{i=1} \norm{y - \y{i}{}{}}^2. \label{eq:13}
\end{multline}
We bound~$\norm{a - \bar{a}}^2$ by plugging
in the definitions of~$a$ and~$\bar{a}$ and 
using the non-expansive property of the projection operator to find 
$
    \norm{a - \bar{a}}^2 
    \leq \sum^N_{j=1} \ \gamma^2_\ell \norm{ \nabla_{\x{}{j}{}} \J{x}{y}{\ell} - \nabla_{\x{}{j}{}} \J{\x{j}{}{}}{\y{j}{}{}}{\ell}}^2.
$
Using the fact that for all~$j \in [N]$ we have
$\big\Vert \nabla_{\x{}{j}{}} \J{x}{y}{\ell} - \nabla_{\x{}{j}{}} \J{\x{j}{}{}}{\y{j}{}{}}{\ell} \big\Vert \leq \big\Vert \nabla_{x} \J{x}{y}{\ell} - \nabla_{x} \J{\x{j}{x}{}}{\y{j}{}{}}{\ell} \big\Vert$
and the Lipschitz property of~$\nabla_{x} \J{\cdot}{\cdot}{\ell}$, we then reach the upper bound 
$
    \norm{a - \bar{a}}^2 \leq \gamma^2_\ell \Lt^2  \sum^N_{j=1} \ \norm{ \big(x,y\big) - \big(\x{j}{}{},\y{j}{}{}\big)}^2. 
$
The result follows by applying this to~\eqref{eq:13}, 
using~$\norm{ \big(x,y\big) - \big(\x{j}{}{},\y{j}{}{}\big)}^2 = \norm{ x - \x{j}{}{}}^2 + \norm{ y - \y{j}{}{}}^2$, and combining like terms. \hfill $\blacksquare$
%

\begin{lemma} \label{lem:qBound}
For all~$\ind{\ell} \in \T$ and for all~$k \geq 0$, we have
    ${\dt{k}{\ell} \leq B^2 m \norm{C}^2  \bt{k}{\ell}}$.
\end{lemma}
\emph{Proof:} 
Suppose agent~$i$ measures its output at time~$k$ and 
suppose that its most recent prior measurement was taken at time~$\mu^i_i(k)$.
Then~$\y{i}{i}{k} =  \y{}{i}{\mu^i_i(k)}$. We know that~${\y{}{i}{k}= C_{i*} x(k)}$, 
and thus 
    $\y{}{i}{k} - \y{}{i}{\mu^i_i(k)} = C_{i*} \x{}{}{k} - C_{i*} \x{}{}{\mu^i_i(k)}$.
Taking the norm squared of both sides and using~$\norm{C_{i*}}^2 \leq \norm{C}^2$ 
leads to the bound
    $\norm{\y{}{i}{k} - \y{}{i}{\mu^i_i(k)}}^2 
                                              \leq \norm{C}^2 \sum_{j=1}^N \ \Big\Vert \sum^{k-1}_{\tau = \mu^i_i(k)} \s{j}{\tau}{\ell} \Big\Vert^2$.
Using the triangle inequality and~$ \Big( \sum^N_{i=1} z_i \Big)^2 \leq N \sum^N_{i=1} z_i^2$  next gives us the bound 
    $\norm{\y{}{i}{k} - \y{}{i}{\mu^i_i(k)}}^2 \leq \norm{C}^2 \sum_{j=1}^N \ \big( k-\mu^i_i(k) \big) \sum^{k-1}_{\tau = \mu^i_i(k)} \norm{\s{j}{\tau}{\ell}}^2$.
By partial asynchrony~$k-\mu^i_i(k) \leq B$ and~$k-B \leq \mu^i_i(k)$, and thus
    $\norm{\y{}{i}{k} - \y{}{i}{\mu^i_i(k)}}^2 \leq \norm{C}^2 B \sum_{j=1}^N \  \sum^{k-1}_{\tau = k-B} \norm{\s{}{\tau}{\ell}}^2$,
where we have used~$\sum^N_{j=1}\norm{\s{j}{\tau}{\ell}}^2 = \norm{\s{}{\tau}{\ell}}^2$. 
By definition of~$\bt{k}{\ell}$ and~$\q{j}{k}$ we can rewrite this as
    $\norm{\q{j}{k}}^2 \leq \norm{C}^2 B \bt{k}{\ell}$.
Taking the sum over the~$m$ outputs on both sides we arrive at
    $\sum^{m}_{j=1} \ \norm{\q{j}{\tau}}^2 = \|q(\tau)\|^2 \leq m \norm{C}^2 B \bt{\tau}{\ell}$.
The result follows by taking the sum from~$\tau = k-B$ to~$\tau = k-1$ on both sides. \hfill $\blacksquare$
\begin{lemma} \label{lem:together}
Take~$\gamma_{\ell} \in (0, 1)$ for all~$\ell\geq 0$.  
For all~$ \ind{\ell} \in\T$ and~$k \in \{\eta_{\ell-1}, \ldots, \eta_\ell - B\}$ we have
\begin{multline}
    \J{\x{}{}{k+B}}{\y{}{}{k+B}}{\ell} - \J{\x{*}{}{\ind{\ell}}}{\y{*}{}{\ind{\ell}}}{\ell} \leq \\
    \big( A_1+A_3+A_6 \big( A_2+A_4 \big) + 1 \big) \sum^{k+B-1}_{\tau = k} \  \norm{\s{}{\tau}{\ell}}^2 \\
    + \big( A_1+A_4+A_5+A_6 \big( A_2+A_4 \big) \big) \sum^{k-1}_{\tau = k-B} \  \norm{\s{}{\tau}{\ell}}^2,
\end{multline}
where we define the constants~$A_1 = N^2M_{\ell}\gamma^2_\ell \Lt^2 + \frac{3}{2} N \Lx^2$ 
and~$A_2 = N^2M_{\ell}\gamma^2_\ell \Lt^2 + \frac{3}{2} N \Ly^2 \norm{C}^2$,
in addition to 
\begin{align}
    A_3 &= NM_{\ell} + N \bigg( \frac{3}{2} K_{\ell}^2 +\frac{5}{2} \bigg) \frac{\lambda^2}{\gamma_{\ell}^2}  \cdot \frac{ 1 + 2\gamma_\ell \big( 1+\gamma_\ell \big) \Lt^2 B N}{1-\gamma_\ell} \\
    A_4 &=  N M_{\ell} + N \bigg( \frac{3}{2} K_{\ell}^2  +\frac{5}{2} \bigg) \frac{\lambda^2}{\gamma_{\ell}^2} \cdot \frac{2\gamma_\ell \big( 1+\gamma_\ell \big) \Lt^2 B N}{1-\gamma_\ell} \\
    A_5 &= N B \frac{\Lx + \Ly N \norm{C}^2}{2} \\
    K_{\ell} &= \Lx + \Ly\|C\|^2, \, M_{\ell} = \frac{3}{2}\left(K_{\ell} + \frac{1}{\gamma_{\ell}}\right)^2 + \frac{5}{2}. 
\end{align}
\end{lemma}
\emph{Proof:}
Fix any~$k \geq 0$. For each~$i$, let~$k^i$ denote the smallest element of~$K^i$ exceeding~$k$. Then we can write the equation
    $\x{i}{i}{k^i+1} = \Pi_{\X_i} \big[ \x{i}{i}{k} - \gamma_\ell \nabla_{\x{}{i}{}} \J{\x{i}{}{k^i}}{\y{i}{}{k^i}}{\ell}  \big]$.
Next we 
apply Lemma~\ref{lem:technical} with~$x=\x{}{}{k}$,~$x^i=\x{i}{}{k^i}$ for ${i \in [N]}$,~$y=\y{}{}{k}$,~$y^i=\y{i}{}{k^i}$ for ${i \in [N]}$, ${\bar{x}=\x{*}{}{\ind{\ell}}}$, ${\bar{y}=\y{*}{}{\ind{\ell}}=C\x{*}{}{\ind{\ell}}}$, 
$a \in \R^n$ with~$a_i = \x{i}{i}{k^i+1}$ for~$i \in [N]$, $b = Ca$
and 
$\bar{a}(k)=\Pi_{\X} \big[ \x{}{}{k} - \gamma_\ell \nabla_{x} \J{\x{}{}{k}}{\y{}{}{k}}{\ell}  \big]$. 
Then 
\begin{multline}
    \J{a}{b}{\ell}-\J{\x{*}{}{\ind{\ell}}}{\y{*}{}{\ind{\ell}}}{\ell} \leq \\ 
    \bigg( N M_{\ell}\gamma^2_\ell \Lt^2 + \frac{3}{2} \Lx^2 \bigg) \sum^N_{i=1} \  \norm{x(k) - \x{i}{}{k^i}}^2 \\
    + \bigg(NM_{\ell}\gamma^2_\ell \Lt^2 + \frac{3}{2} \Ly^2 \norm{C}^2 \bigg) \sum^N_{i=1} \norm{ y(k) - \y{i}{}{k^i}}^2 \\
     + N \bigg( \frac{3}{2} M_{\ell} + \frac{5}{2} \bigg) \norm{\bar{a}(k) - x(k)}^2 
     + N \bigg( \frac{3}{2} K_{\ell}^2 + \frac{5}{2} \bigg) \norm{\x{*}{}{\ind{\ell}} - x(k)}^2.
\end{multline}
Applying Lemma~\ref{lem:errorBound}, using~$0 < \gamma_\ell < 1$, and simplifying gives
\begin{multline} \label{eq:new41}
    \J{a}{b}{\ell}-\J{\x{*}{}{\ind{\ell}}}{\y{*}{}{\ind{\ell}}}{\ell} \leq \\
\bigg( NM_{\ell}\gamma^2_\ell \Lt^2 + \frac{3}{2} \Lx^2 \bigg) \sum^N_{i=1} \  \norm{x(k) - \x{i}{}{k^i}}^2 \\
+ \bigg(NM_{\ell}\gamma^2_\ell \Lt^2 + \frac{3}{2} \Ly^2 \norm{C}^2 \bigg) \sum^N_{i=1} \ \norm{ y(k) - \y{i}{}{k^i}}^2 \\
    + \Bigg( N M_{\ell} 
    + N \bigg( \frac{3}{2} K_{\ell}^2 + \frac{5}{2} \bigg)
    \frac{\lambda^2}{\gamma_{\ell}^2} \Bigg) \norm{\bar{a}(k) - x(k)}^2 .
\end{multline}
We note that~$\|x(k) - x^i(k^i)\|^2 = \sum_{j=1}^{N} \|x_j(k) - x^i_j(k^i)\|^2$ and similar for~$y$. 
Assumption~\ref{ass:partialAsynch} gives us both~$k-B \leq \tau^i_j (k^i) \leq k + B - 1$ and~$k-B \leq \mu^i_j (k^i) \leq k + B - 1$ for all~$i,j \in [N]$. 
Using these facts and the definitions in~\eqref{eq:3} and \eqref{eq:4}, 
we rewrite~\eqref{eq:new41} as 
\begin{multline} \label{eq:17}
    \J{a}{b}{\ell}-\J{\x{*}{}{\ind{\ell}}}{\y{*}{}{\ind{\ell}}}{\ell} \leq 
    A_1 \sum^{k+B-1}_{\tau = k-B} \  \norm{\s{}{\tau}{\ell}}^2 \\
    + A_2 \sum^{k+B-1}_{\tau =k-B} \ \norm{ \q{}{\tau}}^2 
    + NM_{\ell}\Bigg(1 + \frac{\lambda^2}{\gamma_{\ell}^2} \Bigg) \norm{\bar{a}(k) - x(k)}^2.
\end{multline}
We also know that
    $\x{}{i}{k+B} - \x{}{i}{k^i+1} = \sum^{k+B-1}_{\tau = k^i+1} \s{i}{\tau}{\ell}$.
By definition of~$a_i=\x{}{i}{k^i+1}$ we have
\begin{align}
    \x{}{}{k+B} = a + \begin{pmatrix} \sum^{k+B-1}_{\tau = k^1+1} \s{1}{\tau}{\ell} \\ \vdots \\ \sum^{k+B-1}_{\tau = k^N+1} \s{N}{\tau}{\ell}.
    \end{pmatrix}
\end{align}
Then we have
    $\f{\x{}{}{k+B}}{\ell} = \f{a + v}{\ell}$,
where we have used~$v \coloneqq \begin{pmatrix}
        \sum^{k+B-1}_{\tau=k^1 + 1} \s{1}{\tau}{\ell}^T &
        \hdots &
        \sum^{k+B-1}_{\tau=k^N + 1} \s{N}{\tau}{\ell}^T
    \end{pmatrix}^T$. 

From Lemma~\ref{lem:descent}, we find
\begin{multline}
    \f{\x{}{}{k+B}}{\ell} \leq \\ \f{a}{\ell} + \sum^N_{i=1} \sum^{k+B-1}_{\tau = k^i + 1} \s{i}{\tau}{\ell}^T \nabla_{x_i} \f{a}{\ell} + \frac{\Lx}{2} \big\Vert v \big\Vert^2.
\end{multline}
Adding~$\nabla_{\x{}{i}{}} \f{\x{i}{}{\tau}}{\ell} - \nabla_{\x{}{i}{}} \f{\x{i}{}{\tau}}{\ell}$ and rearranging we get
\begin{multline}
    \f{\x{}{}{k+B}}{\ell} \leq \f{a}{\ell} 
    + \sum^N_{i=1} \sum^{k+B-1}_{\tau = k^i + 1} \s{i}{\tau}{\ell}^T \nabla_{\x{}{i}{}} \f{\x{i}{}{\tau}}{\ell} \\
    + \sum^N_{i=1} \sum^{k+B-1}_{\tau = k^i + 1}  \s{i}{\tau}{\ell}^T \Big( \nabla_{x_i} \f{a}{\ell} - \nabla_{\x{}{i}{}} \f{\x{i}{}{\tau}}{\ell} \Big) + \frac{\Lx}{2} \big\Vert v \big\Vert^2.
\end{multline}
Using the Lipschitz property of~$\nabla_{\x{}{}{}} \f{\cdot}{\ell}$ 
\begin{multline}
    \f{\x{}{}{k+B}}{\ell} \leq \f{a}{\ell} 
    + \sum^N_{i=1} \sum^{k+B-1}_{\tau = k^i + 1} \s{i}{\tau}{\ell}^T \nabla_{\x{}{i}{}} \f{\x{i}{}{\tau}}{\ell} \\    
    + \Lx \sum^N_{i=1} \sum^{k+B-1}_{\tau = k^i + 1} \big\Vert \s{i}{\tau}{\ell} \big\Vert  \big\Vert a - \x{i}{}{\tau} \big\Vert 
    + \frac{\Lx}{2} \big\Vert v \big\Vert^2. \label{eq:42}
\end{multline}
By the triangle inequality and
partial asynchrony we next find that  
    $\big\Vert a - \x{i}{}{\tau} \big\Vert \leq \sum^N_{j=1} \sum^{k+B-1}_{\zeta= \tau - B} \big\Vert \s{j}{\zeta}{\ell} \big\Vert$.
Using this in~\eqref{eq:42} gives
\begin{multline}
    \f{\x{}{}{k+B}}{\ell} \leq \f{a}{\ell} 
    + \sum^N_{i=1} \sum^{k+B-1}_{\tau = k^i + 1} \s{i}{\tau}{\ell}^T \nabla_{\x{}{i}{}} \f{\x{i}{}{\tau}}{\ell} \\
    + \Lx \sum^N_{i=1} \sum^{k+B-1}_{\tau = k^i + 1}   \sum^{k+B-1}_{\zeta = \tau - B} \sum^N_{j=1} \big\Vert \s{i}{\tau}{\ell} \big\Vert \big\Vert \s{j}{\zeta}{\ell} \big\Vert
    + \frac{\Lx}{2} \big\Vert v \big\Vert^2.
\end{multline}
Using 
${\norm{\s{i}{\tau}{\ell}} \cdot \norm{\s{j}{\zeta}{\ell}}  \leq  \frac{1}{2} \big( \norm{\s{i}{\tau}{\ell}}^2 + \norm{\s{j}{\zeta}{\ell}}^2} \big)$,
along with~$k \leq k^i + 1$ for all~$i \in [N]$ 
and~$k \leq \tau$, we have
%
\begin{multline}
    \f{\x{}{}{k+B}}{\ell} \leq \f{a}{\ell} + \sum^N_{i=1} \sum^{k+B-1}_{\tau = k^i + 1} \s{i}{\tau}{\ell}^T \nabla_{\x{}{i}{}} \f{\x{i}{}{\tau}}{\ell} \\
    + \frac{\Lx}{2} \Bigg( 2B N \!\!\!\!\! \sum^{k+B-1}_{\tau = k} \!\!\big\Vert\s{}{\tau}{\ell} \big\Vert^2 \!+
    B N \!\!\!\! \sum^{k+B-1}_{\zeta = k - B} \!\! \big\Vert\s{}{\zeta}{\ell}\big\Vert^2 \Bigg) + \frac{\Lx}{2} \big\Vert v \big\Vert^2, \label{eq:fDescent3}
\end{multline}
where we have also used~$\sum^N_{i=1} \norm{\s{i}{k}{\ell}}^2 = \norm{\s{}{k}{\ell}}^2$. 

Next, we continue by bounding~$\g{\y{}{}{k+B}}{\ell}$ in a similar manner. 
Using~$b_j=C_{j*} a$ and~$\y{j}{j}{\mu^i_j(\tau)} = C_{j*} \x{}{}{\mu^i_j(\tau)}$ where~${C_{j*} \in \R^{m_j \times n}}$,
we follow the same steps used to bound~$f\big(x(k+B; t_{\ell}\big)$. Doing this gives
\begin{multline}
    \g{\y{}{}{k+B}}{\ell} \leq \g{b}{\ell} 
    + \sum^N_{i=1} \sum^{k+B-1}_{\tau = k^i +1}  \s{i}{\tau}{\ell}^T C_i^T \nabla_{\y{}{}{}} \g{\y{i}{}{\tau}}{\ell} \\
    + \frac{\Ly N \norm{C}^2}{2} \Bigg( 2B N \sum^{k+B-1}_{\tau = k} \big\Vert\s{}{\tau}{\ell} \big\Vert^2  
    + B N \sum^{k+B-1}_{\zeta = k - B} \big\Vert\s{}{\zeta}{\ell}\big\Vert^2 \Bigg) \\
    + \frac{\Ly \norm{C}^2}{2} \norm{v}^2. \label{eq:gDescent3}
\end{multline}
Adding together Equations~\eqref{eq:fDescent3} and~\eqref{eq:gDescent3} we obtain
\begin{align}
    &\J{\x{}{}{k+B}}{\y{}{}{k+B}}{\ell} - \J{a}{b}{\ell} \leq \\
    &\sum^N_{i=1} \sum^{k+B-1}_{\tau = k^i + 1} \s{i}{\tau}{\ell}^T \nabla_{\x{}{i}{}}\J{\x{i}{}{\tau}}{\y{i}{}{\tau}}{\ell} \\
    & + \frac{2B N \Lx + 2B \Ly N^2 \norm{C}^2}{2}  \sum^{k+B-1}_{\tau = k} \big\Vert\s{}{\tau}{\ell} \big\Vert^2  \\
    &+ \frac{ B N \Lx + B \Ly N^2 \norm{C}^2}{2}  \sum^{k+B-1}_{\zeta = k - B}  \big\Vert\s{}{\zeta}{\ell}\big\Vert^2  \\
    &+ \Bigg( \frac{\Lx + \Ly \norm{C}^2}{2} \Bigg) \big\Vert v \big\Vert^2. \label{eq:47}
\end{align}
By definition of~$v$ we have
    $\big\Vert v \big\Vert^2 = \sum^N_{j=1} \big\Vert \sum^{k+B-1}_{\tau = k^j+1} \s{j}{\tau}{\ell} \big\Vert^2$.
Using the triangle inequality and~$ \Big( \sum^N_{i=1} z_i \Big)^2 \leq N \sum^N_{i=1} z_i^2$ gives 
    $\big\Vert v \big\Vert^2 \leq \sum^N_{j=1}  \Big( \big(k+B -1 \big) - \big( k^j+1 \big) + 1 \Big)  \sum^{k+B-1}_{\tau = k^j+1} \big\Vert \s{j}{\tau}{\ell} \big\Vert^2$.
By partial asynchrony, for each~$j \in [N]$, we have the inequalities $\big(k+B -1 \big) - \big( k^j+1 \big) + 1 \leq B$ and~$k \leq k^j + 1$. Thus, we may write
    $\big\Vert v \big\Vert^2 \leq B \sum^N_{j=1}   \sum^{k+B-1}_{\tau = k} \big\Vert \s{j}{\tau}{\ell} \big\Vert^2$.
Next, using~$\sum^N_{i=1} \norm{\s{i}{k}{\ell}}^2 = \norm{\s{}{k}{\ell}}^2$ gives
    $\big\Vert v \big\Vert^2 \leq B   \sum^{k+B-1}_{\tau = k} \big\Vert \s{}{\tau}{\ell} \big\Vert^2$.
Applying this to~\eqref{eq:47} and using Lemma~\ref{lem:projection} yields
\begin{align}
    &\J{\x{}{}{k+B}}{\y{}{}{k+B}}{\ell} - \J{a}{b}{\ell} \\ 
    &\leq \Bigg( {-\frac{1}{\gamma_{\ell}}} + \frac{3B N \Lx + 3B \Ly N^2 \norm{C}^2}{2} \\
    &\qquad\qquad\qquad+ \frac{B\Lx + B \Ly \norm{C}^2}{2}  \Bigg) \sum^N_{i=1} \sum^{k+B-1}_{\tau = k^i + 1} \big\Vert \s{}{\tau}{\ell} \big\Vert^2 \\
    &  + \frac{ B N \Lx + B \Ly N^2 \norm{C}^2}{2}  \sum^{k-1}_{\zeta = k - B}  \big\Vert\s{}{\zeta}{\ell}\big\Vert^2 + \sum^N_{i=1} \big\Vert \s{i}{k^i}{\ell} \big\Vert^2 \\
    &\leq N B \frac{\Lx + \Ly N \norm{C}^2}{2}  \sum^{k-1}_{\zeta = k - B}  \big\Vert\s{}{\zeta}{\ell}\big\Vert^2 + \sum^N_{i=1} \big\Vert \s{i}{k^i}{\ell} \big\Vert^2, 
\end{align}
which follows from~$\gamma_\ell <  \frac{2}{\big( 3N+1 \big) B \Lx + \big( 3N^2+1 \big) B \norm{C}^2  \Ly}$. 
Adding the last inequality to~\eqref{eq:17} yields 
\begin{multline} 
    \J{\x{}{}{k+B}}{\y{}{}{k+B}}{\ell} - \J{\x{*}{}{\ind{\ell}}}{\y{*}{}{\ind{\ell}}}{\ell} \leq \\
    \bigg( N^2M_{\ell}\gamma^2_\ell \Lt^2 + \frac{3}{2} N \Lx^2 \bigg) \sum^{k+B-1}_{\tau = k-B} \  \norm{\s{}{\tau}{\ell}}^2 + \sum^N_{i=1} \big\Vert \s{i}{k^i}{\ell} \big\Vert^2 \\
    + \bigg( N^2 M_{\ell}\gamma^2_\ell \Lt^2 + \frac{3}{2} N \Ly^2 \norm{C}^2 \bigg) \sum^{k+B-1}_{\tau =k-B} \ \norm{ \q{}{\tau}}^2    \\
    + \Bigg( N M_{\ell} + N M_{\ell}
    \Big( \lambda \max \big\{1, \gamma_\ell^{-1} \big\} \Big)^2 \Bigg) \norm{\bar{a}(k) - x(k)}^2 \\
    + N B \frac{\Lx + \Ly N \norm{C}^2}{2} \sum^{k-1}_{\tau = k - B} \norm{\s{}{\tau}{\ell}}^2. \label{eq:97}
\end{multline}
Since~$k \leq k^i \leq k+B-1$ and~$\sum^N_{i=1} \big\Vert \s{i}{\tau}{\ell} \big\Vert^2 = \big\Vert \s{}{\tau}{\ell} \big\Vert^2$ we have
\begin{equation}
    \sum^N_{i=1} \big\Vert \s{i}{k^i}{\ell} \big\Vert^2 \leq \sum^N_{i=1}  \sum^{k+B-1}_{\tau=k} \big\Vert \s{i}{\tau}{\ell} \big\Vert^2 \leq \sum^{k+B-1}_{\tau=k} \big\Vert \s{}{\tau}{\ell} \big\Vert^2. \label{eq:87}
\end{equation}
Applying Lemma~\ref{lem:update} and~\eqref{eq:87} to~\eqref{eq:97} yields
\begin{multline}
    \J{\x{}{}{k+B}}{\y{}{}{k+B}}{\ell} - \J{\x{*}{}{\ind{\ell}}}{\y{*}{}{\ind{\ell}}}{\ell} \leq \\
    \bigg( N^2  M_{\ell}\gamma^2_\ell \Lt^2 + \frac{3}{2} N \Lx^2 \bigg) \sum^{k+B-1}_{\tau = k-B} \  \norm{\s{}{\tau}{\ell}}^2 \\
    + \bigg( N^2 M_{\ell}\gamma^2_\ell \Lt^2 + \frac{3}{2} N \Ly^2 \norm{C}^2 \bigg) \sum^{k+B-1}_{\tau =k-B} \ \norm{ \q{}{\tau}}^2    \\
    + N M_{\ell}\Bigg(1 + \frac{\lambda^2}{\gamma_{\ell}^2} \Bigg) 
     \frac{ 1 + 2\gamma_\ell \big( 1+\gamma_\ell \big) \Lt^2 B N}{1-\gamma_\ell}  \sum_{\tau=k}^{k+B-1} \Vert \s{}{\tau}{\ell} \Vert^2 \\
     + \epsilon_{\ell} \left( \sum_{\tau=k}^{k+B-1} \Vert \q{}{\tau} \Vert^2 
     + \sum_{\tau=k-B}^{k-1} \Vert \s{}{\tau}{\ell} \Vert^2 
     + \sum_{\tau=k-B}^{k-1} \Vert \q{}{\tau} \Vert^2\right) \\
     + N B \frac{\Lx + \Ly N \norm{C}^2}{2} \sum^{k-1}_{\tau = k - B} \norm{\s{}{\tau}{\ell}}^2 + \sum^{k+B-1}_{\tau=k} \big\Vert \s{}{\tau}{\ell} \big\Vert^2,
\end{multline}
where~$\epsilon_{\ell} = N M_{\ell}\Big(1 + \frac{\lambda^2}{\gamma_{\ell}^2} \Big) 
     \frac{2\gamma_\ell \big( 1+\gamma_\ell \big) \Lt^2 B N}{1-\gamma_\ell} $.
We conclude by applying Lemma~\ref{lem:qBound} and combining like terms. 
\hfill $\blacksquare$
\begin{lemma} \label{lem:moreBound}
For all~$\ind{\ell} \in \T$ and for all~$k \in \{\eta_{\ell-1}, \ldots, \eta_{\ell} - B\}$, we have 
    $\at{k+B}{\ell} \leq \gamma_\ell^{-2}Q_\ell \bt{k+B}{\ell} 
    + \gamma_\ell^{-1} R_\ell \bt{k}{\ell}$,
where we define
    $Q_\ell = Z_1+ Z_3 + A_6 \big( Z_2 + Z_4 \big) + 1$ and
    $R_\ell = Z_1 + Z_4 + A_5+ A_6 \big(Z_2 + Z_4 \big)$,
where~$A_5$ and~$A_6$ are from Lemma~\ref{lem:together}
and
\begin{multline}
    Z_1  \coloneqq \frac{N}{2} \Big(3 \big\Vert C \big\Vert^4 L_\ell^2 \Ly^2 N +6 \big\Vert C \big\Vert^2 L_\ell^2 \Lx \Ly N 
    +3 L_\ell^2 \Lx^2 N \\ +6 \big\Vert C \big\Vert^2 L_\ell^2 \Ly N  +6 L_\ell^2 \Lx N +8 L_\ell^2 N +3 \Lx^2\Big)   
\end{multline}
\begin{multline}
    Z_2 \coloneqq \frac{N}{2}\Big(3 \big\Vert C \big\Vert^4 L_\ell^2 \Ly^2 N +6 \big\Vert C \big\Vert^2 L_\ell^2 \Lx \Ly N +3 \big\Vert C \big\Vert^2 \Ly^2\\
    + 6 \big\Vert C \big\Vert^2 L_\ell^2 \Ly N +3 L_\ell^2 \Lx^2 N +6 L_\ell^2 \Lx N +8 L_\ell^2 N  \Big)
\end{multline}
\begin{multline}
    Z_3 \coloneqq \frac{N}{2} \Big[
    (1+\lambda^2) \Big( 36  B \big\Vert C \big \Vert^4 L_\ell^2 \Ly^2 N 
      + 9 \big\Vert C \big \Vert^4 \Ly^2 \\
      + 72  B \big\Vert C \big \Vert^2 L_\ell^2 \Lx \Ly N 
      + 18 \big\Vert C \big \Vert^2 \Lx \Ly 
      + 36  B L_\ell^2 \Lx^2 N \\
      + 9 \Lx^2 \Big)
      + 72  B \big\Vert C \big \Vert^2 L_\ell^2 \Ly N 
      + 72  B L_\ell^2 \Lx N  
      + 60  B L_\ell^2 N  \lambda^2 \\
      + 96  B L_\ell^2 N +18  \Ly \big\Vert C \big \Vert^2 
      +18 \Lx 
      +15 \lambda^2+24\Big]
\end{multline}
\begin{align}
    &Z_4 \coloneqq N^2 L_\ell^2 B \Big(
    \big( 1 + \lambda^2 \big) \big( 18  \big\Vert C \big\Vert^4 \Ly^2
    + 36  \big\Vert C \big\Vert^2 \Lx \Ly  \\
    &
    + 18  \Lx^2 \big)
    +36  \big\Vert C \big\Vert^2 \Ly 
    +36  \Lx
    +30  \lambda^2
    +48
    \Big).
\end{align}
\end{lemma}
\emph{Proof:} From Lemma~\ref{lem:together}, for any~$k \in \{\eta_{\ell-1}, \ldots, \eta_\ell - B\}$ we have
\begin{multline} 
    \at{k+B}{\ell} \leq \big( A_1+A_3+A_6 \big( A_2+A_4 \big) + 1 \big) \bt{k+B}{\ell} \\
        + \big( A_1+A_4+A_5+A_6 \big( A_2+A_4 \big) \big) \bt{k}{\ell} , \label{eq:91}
\end{multline}
where~$A_1, \dots , A_6 > 0$ are from Lemma~\ref{lem:together}.
Using~$ \gamma_\ell < 1$ gives
\begin{multline}
    A_1 \leq \frac{N}{2} \Big(3 \big\Vert C \big\Vert^4 L_\ell^2 \Ly^2 N +6 \big\Vert C \big\Vert^2 L_\ell^2 \Lx \Ly N +3 \Lx^2\\
    +6 \big\Vert C \big\Vert^2 L_\ell^2 \Ly N +3 L_\ell^2 \Lx^2 N +6 L_\ell^2 \Lx N +8 L_\ell^2 N \Big) . \label{eq:92}
\end{multline} 
Similarly, using~$\gamma_\ell < 1$ lets us upper bound~$A_2$ as 
\begin{multline}
    A_2 \leq \frac{N}{2}\Big(3 \big\Vert C \big\Vert^4 L_\ell^2 \Ly^2 N +6 \big\Vert C \big\Vert^2 L_\ell^2 \Lx \Ly N +3 \big\Vert C \big\Vert^2 \Ly^2\\
    6 \big\Vert C \big\Vert^2 L_\ell^2 \Ly N +3 L_\ell^2 \Lx^2 N +6 L_\ell^2 \Lx N +8 L_\ell^2 N  \Big). \label{eq:93}
\end{multline}

Taking~$\gamma_\ell < \frac{1}{2}$ we have~$\frac{1}{1-\gamma_\ell} < 1+2 \gamma_\ell$. 
Using this inequality and the fact that~$\max \big\{1, \gamma_\ell^{-1} \big\} = \frac{1}{\gamma_\ell}$ yields
\begin{multline}
    A_3 \leq \frac{1}{2 \gamma_\ell^2}\bigg(N \Big(
    \big(1+\lambda^2 \big) \big( 36  B \big\Vert C \big \Vert^4 L_\ell^2 \Ly^2 N 
      + 9 \big\Vert C \big \Vert^4 \Ly^2 \\
      + 72  B \big\Vert C \big \Vert^2 L_\ell^2 \Lx \Ly N 
      + 18 \big\Vert C \big \Vert^2 \Lx \Ly 
      + 36  B L_\ell^2 \Lx^2 N \\
      + 9 \Lx^2 \big)
      + 72  B \big\Vert C \big \Vert^2 L_\ell^2 \Ly N 
      + 72  B L_\ell^2 \Lx N  + 60  B L_\ell^2 N  \lambda^2 \\
      + 96  B L_\ell^2 N +18  \Ly \big\Vert C \big \Vert^2 
      +18 \Lx 
      +15 \lambda^2
      +24\Big)\bigg). \label{eq:94}
\end{multline}
Lastly, we bound~$A_4$. 
Using~$\gamma_\ell < \frac{1}{2}$ we again have~$\frac{1}{1-\gamma_\ell} < 1+2 \gamma_\ell$. Applying this inequality and the fact that~$\max \big\{1, \gamma_\ell^{-1} \big\} = \frac{1}{\gamma_\ell}$ yields
\begin{align}
A_4 &\leq \frac{1}{\gamma_\ell}\bigg(N^2 L_\ell^2 B \Big(
\big( 1 + \lambda^2 \big) \big( 18  \big\Vert C \big\Vert^4 \Ly^2
+ 36  \big\Vert C \big\Vert^2 \Lx \Ly  \\
&
+ 18  \Lx^2 \big)
+36  \big\Vert C \big\Vert^2 \Ly 
+36  \Lx
+30  \lambda^2
+48
\Big)\bigg). \label{eq:95}
\end{align}
From~\eqref{eq:92}-\eqref{eq:95} we have
    $A_1 \leq Z_1$, 
    $A_2 \leq Z_2$, 
    $A_3 \leq Z_3 \gamma_\ell^{-2}$, and 
    $A_4 \leq Z_4 \gamma_\ell^{-1}$.
Using these in~\eqref{eq:91} with~$\gamma_\ell < 1$ gives the result. \hfill $\blacksquare$


\begin{lemma} \label{lem:alphaBound}
For~$\ind{0}$, we have $\at{0}{0} \leq a_0$ and~${\at{B}{0} \leq a_0}$, 
where
    $a_0 \geq L_{J,0} \big( 1 + \norm{C} \big) \diam(\X)$.    
\end{lemma}

\emph{Proof:} 
By definition of~$\at{0}{0}$ we have
\begin{align}
    \at{0}{0} &= \J{\x{}{}{0}}{\y{}{}{0}}{0} - \J{\x{*}{}{\ind{0}}}{\y{*}{}{\ind{0}}}{0}.   \label{eq:29}
\end{align}
Using this and the Lipschitz continuity of~$\J{\cdot}{\cdot}{\ell}$ 
we have
    $\at{0}{0} \leq L_{J,0} \Big( \big\Vert \x{}{}{0} - \x{*}{}{\ind{0}} \big\Vert + \big\Vert \y{}{}{0} - \y{*}{}{\ind{0}} \big\Vert \Big)$.
Then using~$y=Cx$ next gives
    $\at{0}{0} 
    \leq L_{J,0} \big( 1 + \norm{C} \big) \diam(\X)$,
which follows from~$\big\Vert \x{}{}{0} - \x{*}{}{\ind{0}} \big\Vert \leq \diam(\X)$. 
The same argument can be made to bound~$\at{B}{0}$ by replacing~$\at{0}{}$ with~$\at{B}{}$ in~\eqref{eq:29}. Thus,
    $\at{B}{0} \leq L_{J,0} \big( 1 + \norm{C} \big) \diam(\X)$.
Selecting any~$a_0 \geq L_{J,0} \big( 1 + \norm{C} \big) \diam(\X)$ gives
the result. \hfill $\blacksquare$
\begin{lemma}\label{lem:betaBound2}
For all~$\ind{\ell} \in \T$, we have
 $\bt{\eta_{\ell-1}}{\ell}\leq b_\ell$
 and $\bt{\eta_{\ell-1}+B}{\ell} \leq b_\ell$,
where
    $b_\ell \geq B \diam\left(\X\right)^{2}$.
\end{lemma}

\emph{Proof}: By definition of $\bt{\eta_{\ell-1}}{\ell}$ 
in~\eqref{eq:104}
and~$\s{}{\tau}{\ell}$ in~\eqref{eq:81} we have
    $\bt{\eta_{\ell-1}}{\ell}  = \sum_{\tau=\eta_{\ell-1}-B}^{\eta_{\ell-1}-1} \sum_{i=1}^{N} \|x^i_i(\tau+1) - x^i_i(\tau)\|^2$.
We can bound this as
    ${\bt{\eta_{\ell-1}}{\ell} 
    \leq  \sum_{\tau=\eta_{\ell-1}-B}^{\eta_{\ell-1}-1} \sup_{v, w \in X} \|v - w\|^2}$.
Using~$\sup_{v \in X, w \in X} \|v - w\| \leq \diam\left(\X \right)$ we then find
the bound 
     $\bt{\eta_{\ell-1}}{\ell} \leq B \diam\left(\X \right)^2$.
An analogous argument shows that~$\bt{\eta_{\ell-1} + B}{\ell} \leq B \diam\left(\X\right)^2$.

\subsection{Proof of Theorem~\ref{thm:1}} \label{sec:t1proof}
For convenience we define~$\s{}{\tau}{0}=0$ for~$\tau<0$.
From Lemma~\ref{lem:cost} and Lemma~\ref{lem:moreBound} for any~$k$ such that~$0 \leq k \leq \eta_0 - B$,
\begin{align}
    \at{k+B}{0} &\leq \at{k}{0} - \gamma_0^{-1}D_0 \bt{k+B}{0} + E_0 \bt{k}{0} \label{eq:20} \\
    \at{k+B}{0} &\leq \gamma_0^{-2}F_0 \bt{k+B}{0} + \gamma_0^{-1} G_0 \bt{k}{0}. \label{eq:19}
\end{align}
 \newline
Rearranging~\eqref{eq:19} gives 
     $\bt{k+B}{0} \geq \frac{\at{k+B}{0}  - \gamma_0^{-1} G_0 \bt{k}{0}}{\gamma_0^{-2}F_0}$.
Substituting this into~\eqref{eq:20} and rearranging yields
\begin{align}
    &\bigg( 1 + \frac{\gamma_0 D_0}{F_0} \bigg) \at{k+B}{0} \leq \at{k}{0} + \bigg(  \frac{D_0 G_0}{F_0}
    + E_0 \bigg) \bt{k}{0}. \label{eq:22}
\end{align}
Fix~$k \in \{B, \ldots, \eta_0 - B\}$ and substitute~$k-B$ for~$k$ in~\eqref{eq:20} to find 
\begin{align}
    \bt{k}{0} &\leq \gamma_0 \frac{ \at{k-B}{0} -\at{k}{0} + E_0 \bt{k-B}{0}}{D_0}. \label{eq:21}
\end{align}
Applying~\eqref{eq:21} to the right-hand side of~\eqref{eq:22} and rearranging, we get
\begin{align}
     &\at{k+B}{0} \leq \bigg( 1 + \frac{\gamma_0 D_0}{F_0} \bigg)^{-1}  \Bigg( \bigg(1 - \gamma_0 \bigg( \frac{G_0}{F_0} + \frac{E_0}{D_0} \bigg)\bigg) \\
     & \cdot \at{k}{0} 
    + \gamma_0 \bigg( \frac{G_0}{F_0}
    + \frac{E_0}{D_0} \bigg) \Big( \at{k-B}{0} + E_0 \bt{k-B}{0} \Big) \Bigg). \label{eq:25}
\end{align}
For any integer~$d \geq 2$, we iterate~\eqref{eq:20} to find
\begin{multline}
    \at{k+dB}{0} \leq \at{k}{0} - \big( \gamma^{-1}_0 D_0 - E_0 \big) \sum^{d-1}_{j=1} \bt{k+jB}{0} \\
    -\gamma^{-1}_0 D_0 \bt{k+dB}{0} + E_0 \bt{k}{0}.
\end{multline}
Using~$\gamma_0 < D_0 / E_0$ and the fact that~$\bt{\tau}{0} \geq 0$ for all~$\tau > 0$ gives ${\at{k+dB}{0} \leq \at{k}{0} - \big( \gamma^{-1}_0 D_0 - E_0 \big) \bt{k+B}{0} + E_0 \bt{k}{0}}$.
Using the nonnegativity of~$\at{k+dB}{0}$ we can then write $0 \leq \at{k}{0} - \big( \gamma^{-1}_0 D_0 - E_0 \big) \bt{k+B}{0} + E_0 \bt{k}{0}$, and thus
\begin{align}
   \bt{k+B}{0} &\leq \frac{\gamma_0}{ D_0 - \gamma_0 E_0} \big( \at{k}{0} + E_0 \bt{k}{0} \big). \label{eq:26}
\end{align}

Next, we will use~\eqref{eq:25} and~\eqref{eq:26} to show linear convergence of Algorithm~\ref{alg:myAlg} for~$\ind{0}$. 
Select~$a_0 \geq L_{J,0} \big( 1 + \norm{C} \big) \diam(\X)$ and~${b_0 \geq B \diam\left(\X\right)^{2}}$ according to Lemma~\ref{lem:alphaBound} and Lemma~\ref{lem:betaBound2}. Then
         $\at{0}{0} \leq a_0$, 
         $\at{B}{0} \leq a_0$, 
         $\bt{0}{0} \leq b_0$, 
         and
         $\bt{B}{0} \leq b_0$.
Then~\eqref{eq:t1_alpha} and~\eqref{eq:t1_beta} hold for~$r_0=0,1$. Next, we will use induction to show that~\eqref{eq:t1_alpha} and~\eqref{eq:t1_beta} hold for all~$r_0 \in \N_0$. For 
the inductive hypothesis suppose~\eqref{eq:t1_alpha} and~\eqref{eq:t1_beta} hold for all~$r_0$ up to~$d \geq 1$. 
We obtain from~$\gamma_0 < 1 \big/ \big( \frac{G_0}{F_0} + \frac{E_0}{D_0} \big)$ and~\eqref{eq:25} that
\begin{multline} \label{eq:t1_alpha_bound_start}
     \at{dB+B}{0} \leq \frac{1}{1 + \gamma_0 \frac{D_0}{F_0}}  \Bigg( \bigg(1 - \gamma_0 \bigg( \frac{G_0}{F_0} + \frac{E_0}{D_0} \bigg)\bigg) \at{dB}{0} \\
    + \gamma_0 \bigg( \frac{G_0}{F_0}
    + \frac{E_0}{D_0} \bigg) \Big( \at{dB-B}{0} + E_0 \bt{dB-B}{0} \Big) \Bigg).
\end{multline}
For $\gamma_0 \leq \frac{1}{2 c_0}$ we have~$\pt{-1}{0} \leq 1 + 2c_0\gamma_0$.
Using this and the inductive hypothesis, 
we find 
\begin{multline}
     \at{dB+B}{0} \leq \frac{1}{1 + \gamma_0 \frac{D_0}{F_0}}  \Bigg( 1 + 2 \gamma^2_0 c_0 \bigg( \frac{G_0}{F_0} + \frac{E_0}{D_0} \bigg) \\
     + \gamma_0 \bigg( \frac{G_0}{F_0} + \frac{E_0}{D_0} \bigg) \big( 1 + 2c_0 \gamma_0 \big) E_0 \frac{b_0}{a_0} \Bigg) a_0 \pt{d-1}{0}. 
\end{multline}
Then, using~$\frac{b_0}{a_0} = \frac{D_0}{8 E_0 \big( \frac{G_0}{F_0} + \frac{E_0}{D_0} \big) F_0}$,
$\gamma_0 \leq \frac{1}{2 c_0}$, and~$\gamma_0 \leq \frac{D_0}{8 F_0 \Big( \frac{G_0}{F_0} + \frac{E_0}{D_0} \Big) c_0}$ and simplifying gives 
\begin{equation} 
     \at{dB+B}{0}  \leq \bigg( \frac{2 F_0 + \gamma_0 D_0}{2 F_0 + 2 \gamma_0 D_0} \bigg) a_0 \pt{d-1}{0}.
\end{equation}
Using~$1 - \gamma_0 \frac{ D_0}{2F_0 + 2 \gamma_0 D_0} = \frac{2F_0 + \gamma_0 D_0}{2F_0 + 2 \gamma_0 D_0}$ and~$\gamma_0 < 1$,
we reach
\begin{equation} \label{eq:t1_alpha_bound_end}
     \at{dB+B}{0} \leq a_0 \pt{d}{0},
\end{equation}
and this completes the induction on~\eqref{eq:t1_alpha}. 

To complete the this inductive argument we make a similar argument for~\eqref{eq:t1_beta}. From~\eqref{eq:26} 
and the inductive hypothesis we have
the bound~$\bt{dB+B}{0} \leq \frac{\gamma_0 \big( \frac{a_0}{b_0} + E_0 \big)}{ D_0 - \gamma_0 E_0} b_0 \pt{d-1}{0}$.
Here~$\gamma_0 < \gamma_{\max, 0}$ ensures
${\gamma_0 \big(\frac{a_0}{b_0} + E_0\big) \leq (D_0 - \gamma_0 E_0)(1 - \gamma_0 c_0)}$,
and we obtain 
   $\bt{dB+B}{0} \leq \big( 1 - \gamma_0 c_0 \big) b_0 \pt{d-1}{0}$. 
The definition of~$\pt{}{0}$ gives 
   $\bt{dB+B}{0} \leq b_0 \pt{d}{0}$. \label{eq:89}
Thus,~\eqref{eq:t1_alpha} and~\eqref{eq:t1_beta} hold for all~${r_0 \in \N_0}$. 
From Lemma~\ref{lem:qBound} we have
    $\dt{dB+B}{1} \leq d_{0} \pt{d}{0}$,
where~${d_{0} = B^2 m \norm{C}^2  b_{0}}$. 

\subsection{Proof of Theorem~\ref{thm:2}} \label{sec:t2proof}
We first show that for~$\gamma_1 \in \big( 0,\gamma_{\max,1} \big)$ and~$r_1 \in \N_0$ there holds
\begin{align}
    \at{\eta_{0} + r_1 B}{1} & \leq a_1 \pt{r_1-1}{1} \label{eq:59}\\
    \bt{\eta_{0} + r_1 B}{1} & \leq b_1 \pt{r_1-1}{1} \label{eq:60}\\
    \dt{\eta_{0} + r_1 B}{1} & \leq d_1 \pt{r_1-1}{1}. \label{eq:99}
\end{align}
Following the same steps to reach~\eqref{eq:25} we find
\begin{multline}
    \at{k+B}{1} \leq \bigg( 1 + \frac{\gamma_1 D_1}{F_1} \bigg)^{-1}  \Bigg( \bigg(1 - \gamma_1 \bigg( \frac{G_1}{F_1} + \frac{E_1}{D_1} \bigg)\bigg)\at{k}{1}  \\
    + \gamma_1 \bigg( \frac{G_1}{F_1}
    + \frac{E_1}{D_1} \bigg) \Big( \at{k-B}{1} + E_1 \bt{k-B}{1} \Big) \Bigg), \label{eq:56}
\end{multline}
and following the same steps to reach~\eqref{eq:26}, we find 
\begin{align}
   \bt{k+B}{1} &\leq \frac{\gamma_1}{ D_1 - \gamma_1 E_1} \big( \at{k}{1} + E_1 \bt{k}{1} \big). \label{eq:57} 
\end{align}
We will use~\eqref{eq:56} and~\eqref{eq:57} to prove the convergence of Algorithm~\ref{alg:myAlg} for~$t_1$, i.e., show that~\eqref{eq:59},~\eqref{eq:60}, and~\eqref{eq:99} hold for any~$r_1 \in \N_0$. First, we show 
that~\eqref{eq:59} and~\eqref{eq:60} hold for~$r_1 = 0,1$. Specifically, we select~$a_{1} > 0$ and~$b_{1} > 0$ such that
\begin{align}
    \begin{array}{cc}
         \at{\eta_0}{1} \leq a_1 & \at{\eta_0 + B}{1} \leq a_1  \label{eq:78}
    \end{array} \\
    \begin{array}{cc}
    \bt{\eta_0}{1} \leq b_1  & \bt{\eta_0 + B}{1} \leq b_1 . \label{eq:79}
    \end{array}
\end{align}
By definition of~$\at{\eta_0}{1}$ we may write
\begin{multline} 
    \at{\eta_0}{1} =  
    \J{\x{}{}{\eta_0}}{\y{}{}{\eta_0}}{1} 
    - \J{\x{*}{}{\ind{0}}}{\y{*}{}{\ind{0}}}{0} \\
    + \J{\x{*}{}{\ind{0}}}{\y{*}{}{\ind{0}}}{0} 
    +\J{\x{}{}{\eta_0}}{\y{}{}{\eta_0}}{0} \\    
    - \J{\x{}{}{\eta_0}}{\y{}{}{\eta_0}}{0}
    - \J{\x{*}{}{\ind{1}}}{\y{*}{}{\ind{1}}}{1}. \label{eq:30}
\end{multline}
Using the triangle inequality, Assumption~\ref{ass:time2},
and Theorem~\ref{thm:1} gives
\begin{multline}
    \at{\eta_0}{1} \leq 
    \Delta L_t + a_{0} \pt{r_{0} -1}{0} \\
    + \big\vert \J{\x{*}{}{\ind{0}}}{\y{*}{}{\ind{0}}}{0} 
    - \J{\x{*}{}{\ind{1}}}{\y{*}{}{\ind{1}}}{1} \big\vert.
\end{multline}
Adding~$\J{\x{*}{}{\ind{0}}}{\y{*}{}{\ind{0}}}{1} - \J{\x{*}{}{\ind{0}}}{\y{*}{}{\ind{0}}}{1}$, using the triangle inequality,
and again using Assumption~\ref{ass:time2} yields
\begin{multline}
    \at{\eta_0}{1} \leq  
    2 \Delta L_t + a_{0} \pt{r_{0} -1}{0} \\
    + \big\vert \J{\x{*}{}{\ind{0}}}{\y{*}{}{\ind{0}}}{1} 
    - \J{\x{*}{}{\ind{1}}}{\y{*}{}{\ind{1}}}{1} \big\vert. \label{eq:49}
\end{multline}
Using the Lipschitz continuity of~$\J{\cdot}{\cdot}{1}$, the triangle inequality, and~$y = Cx$, 
we find
\begin{multline}
    \big\vert \J{\x{*}{}{\ind{0}}}{\y{*}{}{\ind{0}}}{1} 
    - \J{\x{*}{}{\ind{1}}}{\y{*}{}{\ind{1}}}{1} \big\vert \\
    \leq L_{J,1} (1 + \|C\|)\Big( \big\Vert \x{*}{}{\ind{0}} - \x{*}{}{\ind{1}} \big\Vert \Big) 
    \leq L_{J,1} \sigma_1 \big( 1 + \big\Vert C \big\Vert \big), \label{eq:48}
\end{multline}
where the last inequality uses Assumption~\ref{ass:time2}.
Using~\eqref{eq:48} in~\eqref{eq:49} yields
\begin{align}
    \at{\eta_0}{1} &\leq  
    a_{0} \pt{r_{0} -1}{0} + 2 \Delta L_t  
    + L_{J,1} \sigma_1 \big( 1 + \big\Vert C \big\Vert \big). \label{eq:80}
\end{align}
Next, we bound~$\at{\eta_0 + B}{1}$. By definition we have
\begin{multline}
    \at{\eta_0 + B}{1} 
     = f\bigg(x(\eta_0) + \sum^{\eta_0 + B - 1}_{\tau = \eta_0} \s{}{\tau}{1}; t_1\bigg) \\
     + g\bigg(y(\eta_0) + \sum^{\eta_0 + B - 1}_{\tau=\eta_0} C \s{}{\tau}{1}; t_1\bigg)
    - \J{\x{*}{}{\ind{1}}}{\y{*}{}{\ind{1}}}{1}.
\end{multline} 
Applying Lemma~\ref{lem:descent}, the bound~$ \big\Vert \sum^N_{i=1} z \big\Vert^2 \leq N \sum^N_{i=1} \Vert z \Vert^2$,
the Cauchy-Schwarz inequality, and the triangle inequality, we find
\begin{multline}
    \at{\eta_0 + B}{1} \leq 
    \f{\x{}{}{\eta_0}}{1} \\
    + \sum^{\eta_0 + B - 1}_{\tau = \eta_0} \big\Vert \s{}{\tau}{1} \big\Vert \big\Vert \nabla_x \f{\x{}{}{\eta_0}}{1} \big\Vert 
    + \frac{L_{x,1} B}{2}  \sum^{\eta_0 + B - 1}_{\tau = \eta_0} \big\Vert \s{}{\tau}{1} \big\Vert^2 \\
    +\g{\y{}{}{\eta_0}}{1} 
    + \Vert C \Vert \sum^{\eta_0 + B - 1}_{\tau = \eta_0} \big\Vert \s{}{\tau}{1} \big\Vert \big\Vert \nabla_y \g{\y{}{}{\eta_0}}{1} \big\Vert \\
    + \frac{L_{y,1} B \Vert C \Vert^2}{2} \sum^{\eta_0 + B - 1}_{\tau = \eta_0} \big\Vert \s{}{\tau}{1} \big\Vert^2
    - \J{\x{*}{}{\ind{1}}}{\y{*}{}{\ind{1}}}{1}.
\end{multline} 
Since~$ \bt{\eta_0 + B}{1} = \sum^{\eta_0 + B - 1}_{\tau = \eta_0} \big\Vert \s{}{\tau}{1} \big\Vert^2$, 
we apply Lemma~\ref{lem:betaBound2} to it, then use the boundedness 
of~$\nabla_x \f{\cdot}{1}$ and~$\nabla_y \g{\cdot}{1}$ from \eqref{eq:gradbound} and 
the fact that~$\Vert \s{}{\cdot}{1} \Vert \leq \diam\left( \X \right)$ to find
\begin{multline}
    \at{\eta_0 + B}{1} \leq 
    \at{\eta_0}{1} + \big(  M_{x,1} +  M_{y,1} \Vert C \Vert \big)  B \diam\left( \X \right) \\
    + \big( L_{x,1} + L_{y,1} \Vert C \Vert^2  \big) \frac{ B^2 \diam \left(\X\right)^{2}}{2}.
\end{multline}
Applying the bound in~\eqref{eq:80} we get
\begin{multline}
    \at{\eta_0 + B}{1} \leq a_{0} \pt{r_{0} -1}{0} + 2 \Delta L_t  
    + L_{J,1} \sigma_1 \big( 1 + \big\Vert C \big\Vert \big) \\
    + \big(  M_{x,1} +  M_{y,1} \Vert C \Vert \big)  B \diam\left( \X \right) \\
    + \big( L_{x,1} + L_{y,1} \Vert C \Vert^2  \big) \frac{ B^2 \diam \left(\X\right)^{2}}{2}. \label{eq:84}
\end{multline}

 Next, we modify the right-hand side of~\eqref{eq:84} to design~$\frac{b_1}{a_1}$.  We observe that~$8E_1\left(\frac{G_1}{F_1} + \frac{E_1}{D_1}\right)F_1 \geq 1$. Indeed, 
we have the bound
    $8E_1\left(\frac{G_1}{F_1} + \frac{E_1}{D_1}\right)F_1 = 8E_1G_1 + 8\frac{E_1^2F_1}{D_1} \geq 1$,
which follows by inspection of~$8E_1G_1$. We note for all~$\gamma_1 \in \left( 0 ,  \frac{2}{\big( 1+ B \big) L_{x,1} + \big( 1 +B N  \big) \norm{C}^2 L_{y,1}} \right)$, we have~$D_1 \leq 1$. Then
we have
    $\frac{D_1}{8E_1\left(\frac{G_1}{F_1} + \frac{E_1}{D_1}\right)F_1} \leq 1$ as well.
We can multiply~$\frac{8E_1\left(\frac{G_1}{F_1} + \frac{E_1}{D_1}\right)F_1}{D_1} \geq 1$ with the last term in~\eqref{eq:84} 
to yield a new upper bound to replace~\eqref{eq:84}. 
%
Therefore if we select 
\begin{align} \label{eq:a1_choice}
    a_1 &= a_{0} \pt{r_{0} -1}{0} + 2 \Delta L_t  
    + L_{J,1} \sigma_1 \big( 1 + \big\Vert C \big\Vert \big) \\
    &
    + \big(  M_{x,1} +  M_{y,1} \Vert C \Vert \big)  B \diam\left( \X \right) \\
    &
    + \frac{8 E_1\left(\frac{G_1}{F_1} + \frac{E_1}{D_1}\right)F_1 B^2 \diam \left(\X\right)^{2}}{2 D_1} \big( L_{x,1} + L_{y,1} \Vert C \Vert^2  \big), 
\end{align}
then we satisfy~\eqref{eq:78}.
Furthermore, selecting~$b_1 = B \diam \left( \X \right)^2$ satisfies~\eqref{eq:79} due to Lemma~\ref{lem:betaBound2}. 
From these selections of~$a_1$ and~$b_1$ we have~$\frac{b_1}{a_1} \leq \frac{D_1}{8 E_1 \big( \frac{G_1}{F_1} + \frac{E_1}{D_1} \big) F_1}$.
These choices of~$a_1, b_1$ satisfy~\eqref{eq:59} and~\eqref{eq:60} for~$r_1 = 0,1$. Next, we prove that~\eqref{eq:59} and~\eqref{eq:60} hold for all~$r_1 \in \N_0$ by induction. For 
the inductive hypothesis suppose that~\eqref{eq:59} and~\eqref{eq:60} hold for all~$r_1$ up to~$d\geq 1$. 
Using the same steps used from~\eqref{eq:t1_alpha_bound_start} to~\eqref{eq:t1_alpha_bound_end}, we have
     $\at{\eta_0 +dB+B}{1} \leq a_1 \pt{d}{1}$,
and this completes induction on~\eqref{eq:59}. 

To complete the inductive argument for~$\beta$ we make a similar argument 
in order to reach~\eqref{eq:60}. From~\eqref{eq:57} we have
\begin{equation} \label{eq:t2_beta1}
   \bt{\eta_0 \!+\! dB \!+\! B}{1} \leq \frac{\gamma_1}{ D_1 \!-\! \gamma_1 E_1} \big( \at{\eta_0 + dB}{1} \!+\! E_1 \bt{\eta_0 \!+\! dB}{1} \big).
\end{equation}
Next, we have
   ${\bt{\eta_0 +dB+B}{1} \leq \frac{\gamma_1 \big( \frac{a_1}{b_1} + E_1 \big)}{ D_1 - \gamma_1 E_1} b_1 \pt{d-1}{1}}$
by the inductive hypothesis. 
We then design ${\gamma_1 < \frac{ \frac{a_1}{b_1} + 2E_1 + D_1 c_1 - \sqrt{\big( \frac{a_1}{b_1} + 2E_1 + D_1 c_1\big)^2 - 4 D_1 E_1 c_1}}{2 E_1 c_1}}$ to ensure that~$\gamma_1 \Big( \frac{a_1}{b_1} + E_1 \Big) \leq \big(  D_1 - \gamma_1 E_1 \big) \big( 1 - \gamma_1 c_1 \big)$ so that we obtain 
   $\bt{\eta_0 + dB+B}{1} \leq \big( 1 - \gamma_1 c_1 \big) b_1 \pt{d-1}{1}$.
Therefore 
\begin{align} \label{eq:100}
   \bt{\eta_0 +dB+B}{1} &\leq b_1 \pt{d}{1}. 
\end{align}
Thus,~\eqref{eq:60} holds for all~$r_1 \in \N_0$. 
Then by Lemma~\ref{lem:qBound} we have
    $\dt{\eta_0 + dB+B}{1} \leq B^2 m \norm{C}^2  b_{1} \pt{d}{1}$.
This proves~\eqref{eq:99}. 

The preceding establishes the base case for the next inductive argument
in which we show that if
    $\at{\eta_{\ell-1} + r_\ell B}{\ell}  \leq a_\ell \pt{r_\ell-1}{\ell}$, 
    $\bt{\eta_{\ell-1} + r_\ell B}{\ell}  \leq b_\ell \pt{r_\ell-1}{\ell}$, and
    $\dt{\eta_{\ell-1} + r_\ell B}{\ell}  \leq b_\ell \pt{r_\ell-1}{\ell}$
hold for a fixed~$\ind{\ell}$, then 
$\dt{\eta_{\ell} + r_{\ell+1} B}{\ell + 1} \leq d_{\ell+1} \pt{r_{\ell+1}-1}{\ell + 1}$
and
\begin{align}
    \at{\eta_{\ell} + r_{\ell+1} B}{\ell + 1} & \leq a_{\ell+1} \pt{r_{\ell+1}-1}{\ell + 1} \label{eq:74}\\
    \bt{\eta_{\ell} + r_{\ell+1} B}{\ell + 1} & \leq b_{\ell+1} \pt{r_{\ell+1}-1}{\ell + 1} \label{eq:75}
\end{align}
also hold. This will complete our proof of Theorem~\ref{thm:2}.

From Lemma~\ref{lem:cost} and Lemma~\ref{lem:moreBound} for any~$k \in \{\eta_{\ell}, \ldots, \eta_{\ell+1} - B\}$,
\begin{align}
    \at{k+B }{\ell+1} &\leq \at{k}{\ell+1} - \gamma_{\ell+1} ^{-1}D_{\ell+1} \bt{k+B }{{\ell+1}} + E_{\ell+1} \bt{k}{{\ell+1}} \label{eq:62} \\
    \at{k+B }{\ell+1} &\leq \gamma_{\ell+1}^{-2}F_{\ell+1} \bt{k+B}{\ell+1} + \gamma_{\ell+1}^{-1} G_{\ell+1} \bt{k}{\ell+1}, \label{eq:63}
\end{align} 
by taking~$\gamma_{\ell+1} < \frac{2}{\big( 1+ B \big) L_{x,\ell+1} + \big( 1 +B N  \big) \norm{C}^2 L_{y,\ell+1}}$.
Rearranging~\eqref{eq:63} to lower-bound~$\bt{k+B}{\ell+1}$ we obtain
\begin{align}
     \bt{k+B}{\ell+1} &\geq \frac{\at{k+B}{\ell+1}  - \gamma_{\ell+1}^{-1} G_{\ell+1} \bt{k}{\ell+1}}{\gamma_{\ell+1}^{-2} F_{\ell+1}}.
\end{align}
Applying this to~\eqref{eq:62} and simplifying gives
\begin{multline}
    \bigg( 1 + \frac{\gamma_{\ell+1} D_{\ell+1}}{F_{\ell+1}} \bigg) \at{k+B}{\ell+1} \leq \at{k}{\ell+1} \\
    + \bigg( \frac{D_{\ell+1} G_{\ell+1}}{F_{\ell+1}}
    + E_{\ell+1} \bigg) \bt{k}{\ell+1}. \label{eq:64}
\end{multline}
Fix~$k \in \{\eta_{\ell} + B, \ldots, \eta_{\ell+1} - B\}$ and replace~$k-B$ by~$k$ in~\eqref{eq:62}.
Then 
    $\bt{k}{\ell+1} \leq \frac{\gamma_{\ell+1}}{D_{\ell+1}} \big( \at{k-B}{\ell+1} -\at{k}{\ell+1} + E_{\ell+1} \bt{k-B}{\ell+1} \big)$.
Applying this to~\eqref{eq:64} gives
\begin{align}
    &\at{k+B}{\ell+1} \leq \bigg( 1 + \frac{\gamma_{\ell+1} D_{\ell+1}}{F_{\ell+1}} \bigg)^{-1}  \Bigg( \bigg(1 -\gamma_{\ell+1} \\
    &\cdot\bigg( \frac{G_{\ell+1}}{F_{\ell+1}} + \frac{E_{\ell+1}}{D_{\ell+1}} \bigg)\bigg) \at{k}{\ell+1} 
    + \gamma_{\ell+1} \bigg( \frac{G_{\ell+1}}{F_{\ell+1}}
    + \frac{E_{\ell+1}}{D_{\ell+1}} \bigg) \\
    & \cdot \Big( \at{k-B}{\ell+1} + E_{\ell+1} \bt{k-B}{\ell+1} \Big) \Bigg). \label{eq:66}
\end{align}
For any~$d \geq 2$, we iterate~\eqref{eq:62} to find
\begin{multline}
    \at{k+dB}{\ell+1} \leq \at{k}{\ell+1} - \gamma^{-1}_{\ell+1} D_{\ell+1} \bt{k+dB}{\ell+1} \\ + E_{\ell+1} \bt{k}{\ell+1} 
      - \big( \gamma^{-1}_{\ell+1} D_{\ell+1} - E_{\ell+1} \big) \sum^{d-1}_{j=1} \bt{k+jB}{\ell+1}.       
\end{multline}
Since~$\gamma_{\ell+1} < D_{\ell+1} / E_{\ell+1}$ and~$\bt{\tau}{\ell+1} \geq 0$ for all~$\tau > 0$, 
we obtain
\begin{multline}
    \at{k+dB}{\ell+1} \leq \at{k}{\ell+1} \\ - \big( \gamma^{-1}_{\ell+1} D_{\ell+1} - E_{\ell+1} \big) \bt{k+B}{\ell+1} + E_{\ell+1} \bt{k}{\ell+1}.
\end{multline}
Using~$\at{k+dB}{\ell+1} \geq 0$,
we rearrange to find 
\begin{equation}
   \bt{k+B}{\ell+1} \leq \frac{\gamma_{\ell+1}}{ D_{\ell+1} - \gamma_{\ell+1} E_{\ell+1}} \big( \at{k}{\ell+1} + E_{\ell+1} \bt{k}{\ell+1} \big). \label{eq:67}
\end{equation}

We will use~\eqref{eq:66} and~\eqref{eq:67} to show that~\eqref{eq:74} 
and~\eqref{eq:75} hold for any~$r_{\ell+1} \in \N_0$.
First, we show that~\eqref{eq:74} and~\eqref{eq:75} hold for~$r_{\ell+1} = 0,1$. 
Specifically, we select~$a_{\ell+1} > 0$ and~$b_{\ell+1} > 0$ such that
\begin{equation}
    \at{\eta_{\ell}}{\ell+1} \leq a_{\ell+1}, \qquad \at{\eta_{\ell} + B}{\ell+1} \leq a_{\ell+1} \label{eq:73}
\end{equation}
\begin{equation}
    \bt{\eta_{\ell}}{\ell+1} \leq b_{\ell+1}, \qquad \bt{\eta_{\ell} + B}{\ell+1} \leq b_{\ell+1}. \label{eq:68}
\end{equation}
Using the definition of~$\at{\eta_{\ell}}{\ell+1}$ we 
repeat the steps used from~\eqref{eq:30} to~\eqref{eq:a1_choice} to find
if we select~$a_{\ell+1}$ such that
\begin{multline}
    a_{\ell+1} = a_{\ell} \pt{r_\ell -1}{\ell} + 2 \Delta L_t  
    + L_{J,{\ell+1}} \sigma_{\ell+1} \big( 1 + \big\Vert C \big\Vert \big) \\
    + \big(  M_{x,{\ell+1}} +  M_{y,{\ell+1}} \Vert C \Vert \big)  B \diam\left( \X \right) \\
    + \frac{8E_{\ell+1}\left(\frac{G_{\ell+1}}{F_{\ell+1}} + \frac{E_{\ell+1} }{ D_{\ell+1}}\right)F_{\ell+1} B^2 \diam \left(\X\right)^{2}}{2 D_{\ell+1}} \\
    \cdot \big( L_{x,\ell+1} + L_{y,\ell+1} \Vert C \Vert^2  \big), 
\end{multline}
then~\eqref{eq:73} holds.
Selecting~$b_{\ell+1} = B \diam \left( \X \right)^2$ satisfies~\eqref{eq:68} due to Lemma~\ref{lem:betaBound2}. 
From these selections of~$a_{\ell+1}$ and~$b_{\ell+1}$ we have the bound
    $\frac{b_{\ell+1}}{a_{\ell+1}} 
    \leq \frac{D_{\ell+1}}{8 E_{\ell+1} \big( \frac{G_{\ell+1}}{F_{\ell+1}} + \frac{E_{\ell+1}}{D_{\ell+1}} \big) F_{\ell+1}}$.
These selections of~$a_{\ell+1}, b_{\ell+1}$ satisfy~\eqref{eq:74} and~\eqref{eq:75} for~$r_{\ell+1} = 0,1$. 

Next, we prove that~\eqref{eq:74} and~\eqref{eq:75} hold for all~$r_{\ell+1} \in \N_0$ by induction. For the inductive hypothesis suppose 
that~\eqref{eq:74} and~\eqref{eq:75} hold for all~$r_{\ell+1}$ up to some~$d\geq 1$. 
Using the same steps used from~\eqref{eq:t1_alpha_bound_start} to~\eqref{eq:t1_alpha_bound_end}
we obtain
     $\at{\eta_\ell + dB+B}{\ell+1} \leq a_{\ell+1} \pt{d}{\ell+1}$,
and this completes the induction on~\eqref{eq:74}.

To complete the inductive argument for~$\beta$ we make a similar argument to 
reach~\eqref{eq:75}. 
Following the same steps used to go from~\eqref{eq:t2_beta1} to~\eqref{eq:100}, 
we find 
   $\bt{\eta_\ell + dB+B}{\ell+1} \leq b_{\ell+1} \pt{d}{\ell+1}$.
Therefore,~\eqref{eq:75} holds for~$\ell+1$.
We also have
    $\dt{\eta_\ell + dB+B}{\ell} \leq d_{\ell+1} \pt{d}{\ell+1}$,
from Lemma~\ref{lem:qBound},
where~$d_{\ell+1} \coloneqq B^2 m \norm{C}^2  b_{\ell+1}$.

\end{appendices}

	\bibliographystyle{IEEEtran}{}
	\bibliography{references}

\end{document}